\documentclass[11pt]{article}
\usepackage{amsmath}
\usepackage{graphicx}
\usepackage{amsthm}
\usepackage{verbatim}
\usepackage{multirow}
\usepackage{color}

\newtheorem{thm}{Theorem}[section] 
\newtheorem{lemma}{Lemma}[section] 

\theoremstyle{definition}
\newtheorem{definition}{Definition}
\newtheorem{remark}{Remark}
\newtheorem{example}{Example}
\newcommand{\bed}{\begin{definition}}
\newcommand{\eed}{\end{definition}}
\newcommand{\sgn}{\mathrm{sgn}}
\newcommand{\hamm}{\mathrm{Hamm}}
\newcommand{\cg}{{\cal G}}
\newcommand{\call}{{\cal I}}

\newcommand{\eps}{\epsilon}

\newcommand{\bitem}{\begin{itemize}}
\newcommand{\eitem}{\end{itemize}}

\newcommand{\goto}{\rightarrow}

\newcommand{\margmax}{\mathrm{argmax}}

\newcommand{\beqn}{\begin{equation}}
\newcommand{\eeqn}{\end{equation}}
\newcommand{\balign}{\begin{align}}
\newcommand{\ealign}{\end{align}}

\newcommand{\cG}{{\cal G}}

\newcommand{\that}{\hat{t}}

\usepackage{amssymb}
\begin{document}

\title{Rare and Weak effects in Large-Scale Inference: methods and phase diagrams}
\author{Jiashun Jin and Zheng Tracy Ke}
\date{Carnegie Mellon University and University of Chicago}

\maketitle

\begin{abstract}
Often when we deal with `Big Data',   the true effects we are interested in are {\it Rare
and Weak} (RW). Researchers  measure a large number of features,  hoping to find perhaps only a small fraction of them to be
relevant to the
research in question;  the effect sizes of the relevant features are individually small so  the true effects are not strong enough to stand out for themselves.

Higher Criticism (HC) and Graphlet Screening (GS) are
two classes of methods that are specifically designed for the Rare/Weak settings.
HC was introduced to determine whether there are any relevant effects in
all  the measured features. More recently,
HC was applied to  classification, where it provides a method for selecting useful predictive features
for trained classification rules.  GS  was introduced as a graph-guided multivariate screening procedure,
and was used for variable selection.

We develop a theoretic framework where we use an {\it Asymptotic
Rare and Weak} (ARW)  model simultaneously controlling the size and prevalence of useful/significant features among the useless/null bulk. At the heart of the ARW model is the so-called
{\it phase diagram}, which is a way to visualize clearly the class of ARW settings where the relevant effects are so
rare or weak that desired goals (signal detection, variable selection, etc.)  are simply impossible to achieve. We show that HC and GS have important advantages over better known procedures and achieve the optimal phase diagrams in a variety of ARW settings.

HC and GS  are flexible ideas that adapt easily to many interesting situations.
We review the basics of these ideas and some of the recent extensions, discuss their connections to existing literature,  and suggest some new  applications  of these ideas.
\end{abstract}




\section{Introduction}  \label{sec:intro}
We are often said to be entering the era of `Big Data'.
High-throughput devices measure thousands
or even millions of different features
per single subject on a daily basis.   Such an activity is the driving force of
many areas of science and technology, including a new branch of statistical
practice which Efron \cite{EfronLSI} calls {\it Large-Scale Inference} (LSI).

In many high-throughput data sets,  the relevant effects are {\it Rare and Weak} (RW).
The researchers expect that only a small fraction of these measured features are relevant
for the research in question,  and the effect sizes of the relevant features are individually
small.  The researchers  do not know in advance which features
are relevant and which are not, so they choose to
measure all features within a certain range systematically and
automatically,  hoping to identify  a small fraction of relevant ones
in the future.

Examples include but are not limited to Genome-Wide Association Study (GWAS) and deep sequencing study,
where we are in the so-called ``large-$p$ small-$n$" paradigm, with $p$ being the number of  SNPs and $n$ being the
number of subjects.  As technology on data acquisition advances, we are able to measure increasingly more features per subject.  However,
the number of relevant features do not grow proportionally, so
the relevant effects are sparse;  in addition, since $n$ is usually not as large as we wish,  the effect sizes of
the relevant features (in the summary statistics, e.g., two-sample $t$-tests) are  individually small.

Effect {\it rarity} is a useful hypothesis proposed as early as 1980's by
Box and Meyer \cite{Box}. Later, this hypothesis was found to be
valid in many applications (e.g.,  wavelet image processing \cite{DonoJohn94}, cosmology and astronomy \cite{Starck},
genetics and genomics \cite{Tibs}) and had inspired a long line of researches, where the common theme
is to exploit sparsity (e.g., \cite{DonoJohn94,DonoJohn95}).

However, these works have been largely focused on the regime where the effects are rare but are individually strong (Rare/Strong), with limited attention to the more challenging Rare/Weak regime;
the latter contains many new phenomena which we have not seen in the Rare/Strong regime, to discover which,
we need new methods and new theoretic frameworks.

Call  a relevant feature a signal and an irrelevant feature a noise. In this paper, we investigate two interconnected  LSI problems
in the Rare/Weak regime.
\bitem
\item  {\it Sparse signal detection}. Given two groups of subjects,  a treatment group  and a control group,  each subject is measured on the same set of features. We are interested in deciding whether there is any difference between two groups. In the Rare/Weak setting, the inter-group difference for any single relevant feature is not significant enough, so we have to combine the strengths of these features.
\item {\it Sparse signal recovery}. We have the same setting as above, but the interest is   to separate
relevant features from the overwhelmingly more irrelevant ones. In the literature, this problem is known as {\it variable selection}.
\eitem
Higher Criticism (HC) and Graphlet Screening (GS)
are two recent classes of methods, specifically designed for
Rare/Weak settings, focusing on detecting and recovering sparse signals, respectively.

HC can be viewed as an approach to combining different $P$-values.
It was originally proposed as a method to detect sparse signals
in the presence of uncorrelated noise.
It was then improved to a more sophisticated form called {\it Innovated HC},
 to deal with the case of correlated noise,
where the new ingredient is to exploit graphic structure of the noise terms.
More recently, HC was applied to classification, where it provides a way to
select useful predictive features for trained classification rules.

HC is closely related to the better known methods of controlling False Discovery Rate (FDR)
\cite{BH95}. However,  the target of the FDR-controlling methods is signal recovery
in the Rare/Strong regime,  and the target of HC is signal detection in the
Rare/Weak regime, where signals are so rare and weak that they are inseparable from the noise terms.

GS is a graph-guided multivariate screening procedure,  originally proposed as an approach to variable selection.  GS can be viewed as an extension of the better known method of Univariate
Screening (US) (also called marginal regression or Sure Screening \cite{GJW,FanLv}), where the  innovation   is to provide  a way to overcome the so-called challenge of `signal cancellation' \cite{RW}  by exploiting the graphic structures of the data.  In the simplest case where different features
are uncorrelated, both GS and US reduce to the classical  method of Hard Thresholding \cite{Wassermanbook}.

GS can also be viewed as an approach to evaluating the significance
of several correlated features:  recently,  it was found to  provide a new approach to
feature ranking.  In many LSI problems (e.g.,
multiple testing, classification, spectral clustering), feature ranking plays a pivotal role,
and GS is potentially useful.

GS is very different from better known variable selection methods of $L^0/L^1$-penalizations,
the goal of which is usually to fully recover Rare/Strong signals.
GS focuses on Rare/Weak signals, where full recovery is usually impossible, so we must
develop methods and theory different from those on $L^0/L^1$-penalizations.

Despite that their goals are seemingly very different,
Innovated HC and GS are closely related to each other, and
the main strategies of both methods are to exploit the graphic structures.

We review the basics of HC and GS
and illustrate possible variations and extensions.
We evaluate HC and GS by developing a theoretic framework  using the {\it Asymptotic Rare and Weak} (ARW) model,
 simultaneously controlling the signal prevalence and signal strengths.
We show that  HC and GS offer theoretical optimality in the ARW model, and have advantages over existing methods.

We visualize the ARW model with the phase diagrams.
The phase space refers to
 the two-dimensional space with axes simultaneously quantifying the signal prevalence and signal strengths. It  partitions into several subregions, where the desired inference (signal detection, variable selection, etc.) is distinctly different; because of the partition of the phase space, we call it the phase diagram.

Phase diagram can be viewed as a new criterion for evaluating optimality
which is particularly appropriate for the ARW model.  Given a problem, different methods may also have different phase diagrams, characterizing the subregions where they succeed and where they fail.  When a method partitions the phase diagram in exactly the same way as the optimal methods, we say it achieves the optimal phase diagram.

We show that in a wide variety of settings,  HC and GS achieve
the optimal phase diagrams for signal detection and signal recovery, respectively.       In many of such settings,
other methods (such as FDR-controlling methods \cite{BH95},  $L^0/L^1$-penalization methods \cite{DonohoStark, Tib96})
do not achieve the optimal phase diagrams.

Note that, however,  HC and GS are flexible ideas and can be applied to many interesting settings. They are  not tied to the ARW model and their advantages over existing methods remain in  much broader settings.

\subsection{Roadmap and highlights}
In Section \ref{sec:ARW}, we introduce the ARW model and the watershed phenomena associated with the problems of   sparse signal detection and  sparse signal recovery, respectively, and visualize them with the phase diagrams.

Section \ref{sec:ARW} does not address the achievability:  which methods (presumably easy-to-use) achieve the optimal phase diagrams. The achievability is addressed in Sections  \ref{sec:white}-\ref{sec:color}, focusing on the cases of uncorrelated noise and correlated noise, respectively.
In Section \ref{sec:white}, we first review the basics of the HC  and show that it achieves the optimal phase diagram
for  signal detection. We then show that the well-known method of Hard Thresholding
achieves the optimal phase diagram for  signal recovery.
We also review the recent literature on the idea of HC.
Section \ref{sec:color} discusses  the case of correlated noise.
For signal detection, we develop HC into the more sophisticated Innovated HC,
and for signal recovery, we use GS.
Phase diagrams and optimality of HC and GS are also investigated.

In Section \ref{sec:future}, we suggest some new applications of HC and GS, supported by some preliminary numeric studies.
In Section \ref{sec:class}, we extend HC as a feature selection method in the context of classification.
We address the classification phase diagram as well as the optimality of HC.

The development of HC and GS exposes several noteworthy ideas; we provide
a road map  to highlight these ideas, with details later.
\begin{itemize}
\item Innovated Transformation (IT) is the key to many methods
(e.g., Innovated HC, GS, and the classification rule in Fan {\it et al}. (2013) \cite{FJY})
that incorporate correlation structures of the noise terms for inference.
Compared to Whitened Transformation (WT) which attempts to
create uncorrelated noise terms by transformations,  IT---though counterintuitive---
yields  larger (post-transformation) Signal-to-Noise Ratios (SNR) than WT,  simultaneously at all  (pre-transformation) signal sites, and so it is preferred.  See Section \ref{subsec:IHC} for details.
\item While it is expected that the optimal phase diagrams critically depend on the local graphic structures,  it  is not  the case, and most parts of the phase diagrams remain the same across a range of very different local graphic structures.   On the other hand,   for  optimal procedures (e.g., GS), we must exploit  local graphic structures.
The well-known $L^0/L^1$-penalization methods do not adequately exploit local graphic structures,   so they do not
achieve the optimal phase diagram, even in very simple settings and even when the tuning parameters are ideally set.
See Sections \ref{subsec:GS}-\ref{subsec:example} for details.
\item  In classification, a prevailing idea is to select a few important predictive features so that the (feature) FDR \cite{BH95} is small. Recent studies reveal something very different: in some Rare/Weak settings,  we must select features in a way so that the FDR is very high, so that we are able to include almost all
useful features for classification. See Section \ref{sec:class}.
\end{itemize}
Other small but noteworthy items include Lemma \ref{colornumber} and Remarks 4, 12.

\section{The ARW model and phase diagrams}
\label{sec:ARW}
In this section, we introduce the ARW model.  We discuss
the watershed phenomena associated with signal detection
and signal recovery, respectively, and visualize them with the
phase diagrams. Discussions on the achievability are long and
deferred to Sections \ref{sec:white}-\ref{sec:color};  we explain
our plan for the discussions on the achievability in Section \ref{subsec:plan}.

The ARW model was first proposed by Donoho and Jin (2004) \cite{DJ04}
for sparse signal detection. More recently,  it was extended to more complicated forms \cite{HJ09} and to different settings including classification \cite{DJ09, IPT, JinPNAS}, variable selection \cite{JiJin, JZZ}, and spectral clustering \cite{JW}; see also \cite{Castro2, Meinshausen, Lin}.

In this paper, we focus on a version of the ARW model that is simple for presentation, yet contains all important ingredients associated with the major insights exposed in the above references.

In such a spirit, we consider a Stein's $p$-normal means model
\begin{equation} \label{Stein1}
Y  =  \beta + z, \qquad z \sim N(0, \Sigma),
\end{equation}
where $\Sigma = \Sigma_{p, p}$ is the covariance matrix.
Denote the precision matrix by
\[
\Omega = \Omega_{p,p} = \Sigma_{p,p}^{-1}.
\]
For simplicity, we usually drop subscripts `$p,p$'.
We assume $\Omega$ is sparse in the (strict) sense that each row of $\Omega$ has relatively few nonzeros.
Such an assumption is only for simplicity in presentation; see \cite{HJ09,JZZ} for discussions on more general $\Omega$.
Model (\ref{Stein1}) may arise from many applications, including the following.
\begin{itemize}
\item {\it Two-group study}.  In the aforementioned two-group study,  $Y_j$, $1 \leq j \leq p$,  can be viewed as the
two-sample $t$-statistic associated with the $j$-th feature.  In many such studies (e.g., Genetic Regulatory Network (GRN)), the precision matrix is sparse \cite{GRN}.
\item {\it Linear models with random designs}. Given $W \sim  N(X \tilde{\beta}, I_n)$,  where the rows of $X$
 are iid samples from the $p$-dimensional distribution $N(0,  \Omega)$, where $\Omega$ is sparse. Such settings can be found in Compressive Sensing \cite{Donoho, DMM09} or Computer Security \cite{Nissim, Fienberg}, where $\Omega = I_p$.
Letting $\widetilde{W} = (1/\sqrt{n}) X'W$ and $\beta = \sqrt{n} \tilde{\beta}$, then approximately $\widetilde{W}  \sim N(\Omega \beta, \Omega)$, which is equivalent to model (\ref{Stein1}); the connection is solidified in  Jin {\it et al}.  (2012) \cite{JZZ}.
\end{itemize}
We assume $\Omega$ is known, as our primary goal is to investigate how the graphic structures of $\Omega$ affect the constructions of the optimal methods and optimal phase diagrams.
Such an assumption is valid in many applications. For example, in the above linear model, $\Omega$ plays the role of
the Gram matrix which is known to us.  When $\Omega$ is unknown but is sparse, it can be estimated by many recent algorithms, such as the glasso \cite{glasso}.
The gained insight here is readily extendable to the case where $\Omega$ is unknown
but can be estimated reasonably well, as only large entries of $\Omega$ have
a major influence on the problems we are interested in.

We choose a somewhat unconventional normalization such that
\begin{equation} \label{diagOmega}
\Omega(i,i) = 1, \qquad 1 \leq i \leq p.
\end{equation}
Denote $d_p^*$ by the maximum number of nonzeros in the rows of $\Omega$:
\[
d_p^* = d_p^*(\Omega) = \max_{1\leq i\leq p} \bigl\{ \# \{1\leq j\leq p: \Omega(i,j)\neq 0\} \bigr\}.
\]
We assume
$d_p^*(\Omega)$ {\it grows slowly enough}:
\begin{equation} \label{degreeOmega}
d_p^*(\Omega)  p^{-\delta} \goto 0,  \qquad \mbox{for any fixed $\delta > 0$}.
\end{equation}

At the same time, fixing $\eps \in (0,1)$ and $\tau > 0$, we  model the vector  $\beta$ by
\begin{equation} \label{twopt}
\beta_i   \stackrel{iid}{\sim}    (1-\eps)  \nu_0+\eps  \nu_{\tau},  \qquad 1 \leq i \leq p,
\end{equation}
where  $\nu_{a}$ is the point mass at $a$.
We are primarily interested in the case where $\eps$ is small and $\tau$ is small or moderately large, so that the signals  (i.e., nonzero entries of $\beta$) are  Rare and Weak.

In our asymptotic framework,  we let $p$ be the driving asymptotic parameter, and  tie $(\eps, \tau)$ to $p$ through fixed parameters $\vartheta$ and $r$:
\begin{equation} \label{Defineeps}
\eps = \eps_p = p^{-\vartheta}, \qquad 0 < \vartheta < 1,
\end{equation}
\begin{equation} \label{Definetau}
\tau  = \tau_p = \sqrt{2 r \log(p)}, \qquad r > 0.
\end{equation}
Note that as $p$ grows, the signals become increasingly sparser.
To counter this effect, we have to let the signal strength parameter $\tau$ grows to $\infty$
slowly, so that the testing problem is non-trivial.
\bed
We call model (\ref{Stein1})--(\ref{Definetau}) the Asymptotic Rare and Weak model $ARW(\vartheta, r, \Omega)$.
\eed
\noindent
See \cite{GJW, HJ08, HJ09, JiJin, JZZ, Ke} for the ARW model in more general forms. The ARW model is subtle even in the case where $\Omega = I_p$; see \cite{DJ04} for example.

\subsection{Detecting Rare and Weak signals}
\label{subsec:detectLB}
We formulate the sparse signal detection problem as a hypothesis testing problem,
where we test
a {\it joint null hypothesis} that  all $\beta_i$'s are $0$:
\begin{equation} \label{testnull}
H_0^{(p)}: \qquad   \beta = 0,
\end{equation}
against a specific complement of the joint null:
\begin{equation} \label{testalt}
H_1^{(p)}: \qquad \mbox{$\beta$ satisfies a Rare and Weak model \eqref{twopt}-\eqref{Definetau}}.
\end{equation}
Note that $Y_i \stackrel{iid}{\sim} N(0,1)$ and $Y_i \stackrel{iid}{\sim} (1 - \eps_p) N(0,1)
+ \eps_p  N(\tau_p, 1)$ under $H_0^{(p)}$ and $H_1^{(p)}$, respectively, so the testing problem (\ref{testnull})-(\ref{testalt}) is also the problem of detecting Gaussian mixtures \cite{DJ04}.

It turns out that there is a watershed phenomenon, meaning that in the two-dimensional phase space $\{(\vartheta, r): 0 < \vartheta < 1, r > 0\}$, there is a curve that partitions the whole phase space into two regions
where the testing problem (\ref{testnull})-(\ref{testalt}) is distinctly different.  Presumably,  the curve should depend on the off-diagonals of $\Omega$ in a complicated way.
Somewhat surprisingly, this is not the case, and the off-diagonals of $\Omega$ do not have a major influence on the
partition.

In detail,  define the {\it standard phase function for detection}
\begin{equation}  \label{phasefunction}
\rho^*(\vartheta) =
\left\{
\begin{array}{ll}
 0,  & 0 <  \vartheta \leq  1/2,  \\
\vartheta - 1/2,   &1/2 <    \vartheta  \leq  3/4,  \\
  (1-\sqrt{1-\vartheta})^2, &3/4 < \vartheta  < 1.
  \end{array}
\right.
\end{equation}
We have the following theorem.
\begin{thm}  \label{thm:Detect}
Fixing $\vartheta \in (0,1)$ and $r > 0$,  consider a sequence of $ARW(\vartheta, r, \Omega)$ indexed by $p$.  If $r < \rho^*(\vartheta)$, then for any sequence of tests that test $H_1^{(p)}$ against $H_0^{(p)}$, the sum of Type I and Type II errors tends to $1$ as $p \goto \infty$.  If $r > \rho^*(\vartheta)$, then there is a test the sum of Type I and Type II errors of which  tends to $0$ as $p \goto \infty$.
\end{thm}
In the simplest case of  $\Omega = I_p$, Theorem \ref{thm:Detect} was proved
in Donoho and Jin (2004) \cite{DJ04}, and also in \cite{Ingster97, Ingster99}.  For general cases $\Omega \neq I_p$,
the second claim follows from Theorem \ref{thm:IHC} on Innovated HC below, and
the proof of the first claim is similar to that in \cite[Theorem 1.1]{FJY} so we skip it;  at the heart of the proof is subtle analysis of the Hellinger distance associated with the testing problem, as well as the following lemma.
\begin{lemma} \label{colornumber}
{\it (Chromatic Number)}.
Fixing $p$ and $K$ such that $1 \leq K < p$, consider a $p \times p$ matrix $\Omega$. If each row of $\Omega$ has no more than $K$ nonzeros, then we can color indices $1, 2, \ldots, p$ in no more than $K$ different colors, so that
for any pair of indices $i,j$ with the same color, $\Omega(i,j) = 0$.
\end{lemma}

\begin{remark}
In (\ref{phasefunction}),  $\rho^*(\vartheta) = 0$ when $0 < \vartheta < 1/2$. This does not mean
that two hypotheses $H_0^{(p)}$ and $H_1^{(p)}$ are asymptotically separable for any $(\eps_p, \tau_p)$; it only means that two hypotheses  can be asymptotically separated  even when $\tau_p \ll \sqrt{\log(p)}$.  See \cite{DJ04} for more discussions.
\end{remark}

\subsection{Recovering Rare and Weak signals}
\label{subsec:recoverLB}
Consider again the ARW model  where
\[
Y = \beta + z, \qquad z \sim N(0, \Sigma),
\]
and $(\beta, \Omega)$ satisfy (\ref{diagOmega})-(\ref{Definetau}).
Now,  the main interest is to separate the nonzero entries of $\beta$ from the zero ones (i.e., signal recovery or variable selection).
For any estimator $\hat{\beta}$, we measure the
errors by the Hamming distance:
\[
h_p(\hat{\beta}, \beta) = E  \sum_{i = 1}^p 1\{ \sgn(\hat{\beta}_i) \neq \sgn(\beta_i)\},
\]
where $\sgn(u)$ denotes the sign of $u$ taking values $0,\pm 1$
and the expectation is taken over the randomness of $\beta$ and $Y$.
Hamming distance is approximately the expected sum of the number of
signals that have been classified as noise and the number of noise that have been classified as signals.
The minimax Hamming distance is then
\begin{equation} \label{eqhamm}
\hamm_p^*(\vartheta, r; \Omega) = \inf_{\hat{\beta}}  \{ h_p(\hat{\beta}, \beta) \}.
\end{equation}

The following short-hand term is frequently used in the paper.
\bed
$L_p > 0$ denotes a generic multi-$\log(p)$ term which may change from occurrence to occurrence and satisfies that  $L_p p^{\delta} \goto \infty$ and $L_p p^{- \delta} \goto 0$ for any $\delta > 0$, as $p \goto \infty$.
\eed
\noindent
The following theorem is proved by Ji and Jin (2011) \cite{JiJin} (see also \cite{JZZ,Ke}).
\begin{thm} \label{thm:recover}
{\it (Lower bound}). Fixing $\vartheta \in (0,1)$ and $r > 0$, consider a sequence of $ARW(\vartheta, r, \Omega)$ indexed by $p$. As $p \goto \infty$,
\begin{equation} \label{eqLB}
\hamm_p^*(\vartheta, r; \Omega)  \left\{
\begin{array}{ll}
\geq L_p \cdot p^{1 - (\vartheta + r)^2/(4r)}, &\qquad r > \vartheta, \\
\sim p^{1 - \vartheta}, &\qquad 0 < r < \vartheta.
\end{array}
\right.
\end{equation}
\end{thm}

Similarly, there is a watershed phenomenon associated with the problem of signal recovery.
This phenomenon can be most conveniently described in the case $\Omega = I_p$, but it holds
much more broadly.

When $\Omega = I_p$,  it is shown in Section \ref{sec:white} below that the lower bound in (\ref{eqLB}) is tight and can be achieved  by Hard Thresholding:
\begin{equation} \label{eqLB2}
\hamm_p^*(\vartheta, r; I_p)      \left\{
\begin{array}{ll}
 =  L_p \cdot p^{1 - (\vartheta + r)^2/(4r)}, &\qquad r > \vartheta, \\
\sim p^{1 - \vartheta}, &\qquad 0 < r < \vartheta.
\end{array}
\right.
\end{equation}
Define the {\it standard phase function for exact recovery}:
\begin{equation} \label{exactboundary}
\rho^*_{exact}(\vartheta; I_p) = (1+ \sqrt{1-\vartheta})^2, \qquad 0<\vartheta <1.
\end{equation}
Let
\[
S_p(\beta) = \{1\leq i \leq p: \beta_i \neq 0\}, \qquad s_p(\beta) = |S_p(\beta)|.
\]
In the ARW model,  with overwhelming probability,
\begin{equation} \label{definesp}
s_p(\beta)   \sim  p \eps_p  =  p^{1 - \vartheta}.
\end{equation}
For any fixed $(\vartheta, r)$, by (\ref{exactboundary})-(\ref{definesp})  and basic algebra, it follows that
\begin{itemize}
\item when $r > \rho^*_{exact}(\vartheta; I_p)$, $\hamm_p^*(\vartheta, r; \Omega)  = o(1)$,
and it is possible to fully recover $S(\beta)$---the support of $\beta$---with overwhelming probabilities.
\item when $\vartheta < r < \rho^*_{exact}(\vartheta; I_p)$,
$1 \ll  \hamm_p^*(\vartheta, r; I_p) \ll p^{1 - \vartheta}$; it is possible to recover most of the signals,
but it is impossible to fully recover $S(\beta)$.
\item  when $0 < r < \vartheta$, $\hamm_p^*(\vartheta, r; I_p)  \sim p^{1 - \vartheta}$,
and signals and noise are merely inseparable.
\end{itemize}

The case $\Omega \neq I_p$ is the same, except for one difference, that we have to replace the function $\rho_{exact}^*(\vartheta; I_p)$ by a more general function
\[
\rho_{exact}^*(\vartheta; \Omega).
\]
We address this case in Section \ref{sec:color}, where
we attack the problem by Graphlet Screening (GS). We derive the optimal rate of convergence of the minimax
Hamming distance $\hamm_p^*(\vartheta, r; \Omega)$, and show that GS is asymptotically  minimax.

In principle,  $\hamm_p^*(\vartheta, r; \Omega)$ may depend on $\Omega$ in a complicate way.
Still, for many sequences of $\Omega = \Omega_{p,p}$ (with careful calibrations, possibly),  there is a constant $c = c(\vartheta, r; \Omega)$ depending on $(\vartheta, r)$ and the calibrations we choose for $\Omega$ such that
\[
\hamm_p^*(\vartheta, r; \Omega) = L_p p^{1 - c(\vartheta, r;  \Omega)},
\]
and equating $c(\vartheta, r; \Omega) = 1$ gives $r=\rho_{exact}^*(\vartheta; \Omega)$.
Examples of such sequences include that of $\Omega = I_p$ for all $p$, a diagonal block-wise
example to be discussed in Section \ref{subsec:example}, and
the long-memory time series model and the change point model discussed in Ke {\it et al}   (2012) \cite{Ke}  as well as in Remark 10.
 In all these cases, we have explicit formulas
for $\rho_{exact}^*(\vartheta; \Omega)$. See details therein.

\subsection{Phase diagrams} \label{subsec:phase}
\begin{figure}[tb]
\centering
\includegraphics[width = 2.8 in]{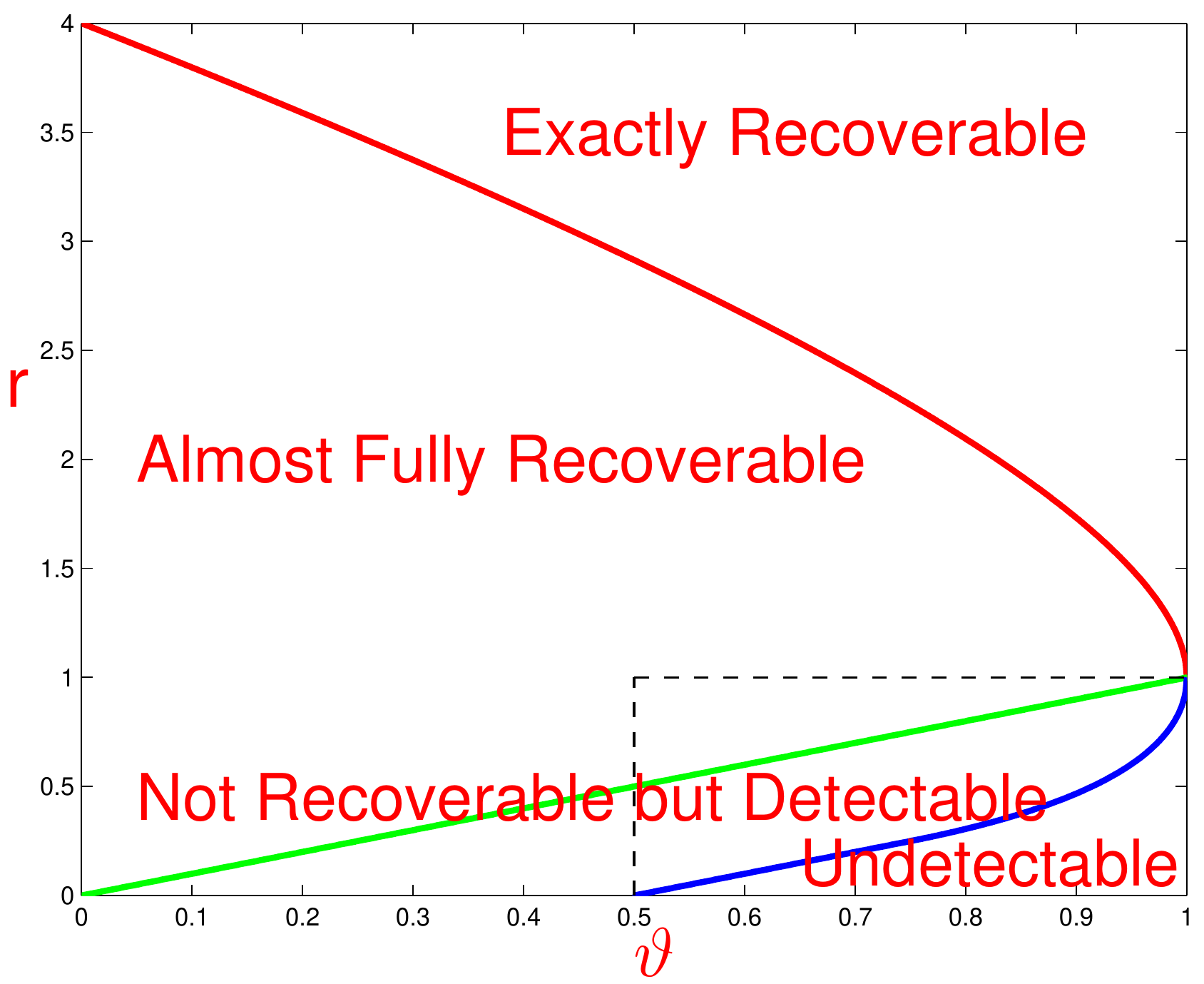}
\includegraphics[width = 2.8 in]{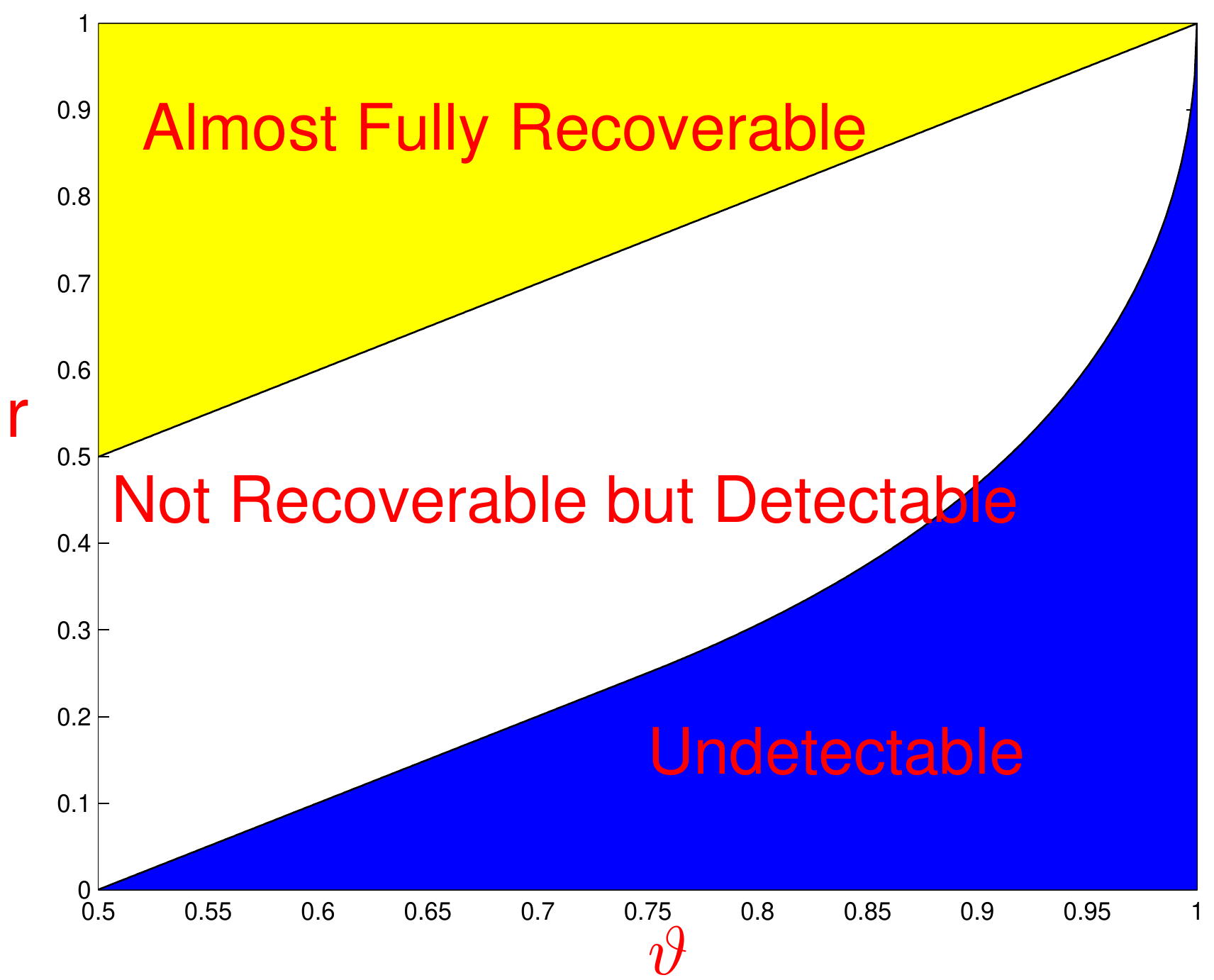}
\caption{Phase diagrams ($\Omega = I_p$).  Left: curves in red, green,  and blue are $r = \rho_{exact}^*(\vartheta; I_p)$, $r = \vartheta$, and $r = \rho^*(\vartheta)$, correspondingly.  Right: enlargement of the region bounded by the dashed lines in the left panel.}
\label{fig:4Phases}
\end{figure}

The preceding results can be  visualized with the phase diagrams.
The two-dimensional phase space $\{(\vartheta, r):  0 < \vartheta < 1, r > 0\}$ is partitioned by
the   three curves
\[
r = \rho^*(\vartheta),  \qquad r = \vartheta, \qquad r = \rho_{exact}^*(\vartheta; \Omega)
\]
into four different subregions, where the inference is distinctly different.
\begin{itemize}
\item {\it Region of Undetectable}: $\{(\vartheta, r):  0 < \vartheta < 1, r <  \rho^*(\vartheta)\}$.
The signals are so rare and weak that it is impossible to detect their existence: the two hypotheses $H_1^{(p)}$ and    $H_0^{(p)}$ are merely inseparable; for any test,   the sum of Type I and Type II errors tends to $1$ as $p \goto \infty$.
\item {\it Region of Not Recoverable but Detectable}:
$\{(\vartheta, r): 0 < \vartheta < 1,  \rho^*(\vartheta) < r < \vartheta\}$. It is possible to have  a test such that the  sum of Type I and Type II errors tends to $0$ as $p \goto \infty$.
However, it is impossible to separate the signals from the noise: as $\hamm_p^*(\vartheta, r; \Omega) \gtrsim p^{1 - \vartheta}$,  the expected Hamming distance of  any estimator is comparable to the total number of signals.
\item {\it Region of Almost Fully Recoverable}: $\{(\vartheta, r): 0 < \vartheta < 1,  \vartheta < r <
\rho_{exact}^*(\vartheta; \Omega)\}$.  It is possible to recover almost all signals but not all of them; the Hamming distance
is much smaller than $p^{1 - \theta}$, but is also much larger than $1$.
\item {\it Region of Exactly Recoverable}: $\{(\vartheta, r):  0 < \vartheta < 1, r > \rho_{exact}^*(\vartheta; \Omega)\}$.  The signals are sufficiently strong so that $\hamm_p^*(\vartheta, r; \Omega)  = o(1)$, and it is possible to have exact recovery with overwhelming probabilities.
\end{itemize}
In the last two sub-regions, it is clearly detectable;  we call the following  region
\begin{equation} \label{Detectable}
\{(\vartheta, r):  0 < \vartheta < 1,    r > \rho^*(\vartheta)\}
\end{equation}
the {\it Region of Detectable}.  Also, note that among these four subregions, only the last two depend  on $\Omega$.
In Figure~\ref{fig:4Phases}, we plot the phase diagrams in the case of $\Omega = I_p$. Phase diagrams for more general cases are discussed in Section \ref{sec:color}.

Phase diagram can be viewed as a new criterion for measuring performances that is particularly appropriate
for Rare/Weak settings. Phase diagram is a flexible idea, which has been extended recently to many different settings, including  large-scale multiple testing \cite{Castro2,HJ09,Wellner2004}, variable selection \cite{JiJin, JZZ, Ke},  classification \cite{DJ09, FJY, IPT, JinPNAS}, spectral clustering and low-rank matrix recovery \cite{JW}, and computer privacy and security \cite{Fienberg}.  See sections below for more discussions.

\subsection{Achievability of the phase diagrams}
\label{subsec:plan}
In the preceding section, we have only said that the optimal phase diagrams are achievable, without referencing
to any specific methods.  It is of primary interest to develop methods---preferably easy-to-implement and is not tied to the ARW model---to achieve the optimal phase diagrams:
\begin{itemize}
\item We say a testing procedure achieves the {\it optimal phase diagram for detection}  if for any $(\vartheta, r)$ in the interior of
Region of Detectable, the power of the procedure tends to $1$ as $p \goto \infty$.
\item We say a variable selection procedure $\hat{\beta}$ achieves the {\it optimal phase diagram for recovery} if $h_p(\hat{\beta}, \beta) \leq  L_p  \cdot  \hamm_p^*(\vartheta, r; \Omega)$  for sufficiently large $p$,
where $L_p$ is the generic multi-$\log(p)$ term as in Definition 2.
\end{itemize}
If a testing procedure and a variable selection procedure  achieve the optimal phase diagrams for detection and recovery, respectively, then they
partition the phase space in exactly the same way as in Section \ref{subsec:phase}.

In Section \ref{sec:white}, we address the achievability for the case $\Omega = I_p$, and show  that
Orthodox Higher Criticism (OHC) and Hard Thresholding achieve the optimal phase diagrams for detection and recovery, respectively. In Section \ref{sec:color}, we address the acheivabilty for more genera $\Omega$, and show that
Innovated HC and GS achieve the optimal phase diagrams for detection and recovery, respectively.
Combining these with Theorems \ref{thm:Detect}-\ref{thm:recover} gives the phase diagrams presented in Section
\ref{subsec:phase}.

\section{Detecting and recovering signals in white noise}
\label{sec:white}
We revisit the problems of signal detection and signal recovery, and show that when $\Omega = I_p$,
HC and Hard Thresholding achieve the optimal phase diagrams for detection and
 recovery, respectively.  We also review the recent applications and extensions of the HC idea. Discussions on general $\Omega$ is in Section \ref{sec:color}.

\subsection{Higher Criticism basics and optimal signal detection (white noise)}  \label{subsec:HCbasics}
Higher Criticism (HC) is a notion that goes back to John Tukey \cite{Tukey76,Tukey89,Tukey94}.  Donoho and Jin (2004) \cite{DJ04} developed HC into a method of combining $P$-values,  and used it to resolve the testing problem (\ref{testnull})-(\ref{testalt}). HC consists of three steps.
\begin{itemize}
\item For each $1 \leq i \leq p$, obtain a $P$-value by $\pi_i = P(N(0,1) \geq Y_i)$.
\item Sort  the $P$-values in the ascending order:
$\pi_{(1)} < \pi_{(2)} < \ldots < \pi_{(p)}$.
\item  The Higher Criticism statistic is then
\begin{equation} \label{OHC}
HC_p^*   =  \max_{\{ 1 \leq i \leq p/2 \}}  HC_{p, i}, \qquad   \mbox{where}  \;\;    HC_{p, i} \equiv  \sqrt{p}  \frac{(i/p)  - \pi_{(i)}}{\sqrt{\pi_{(i)} (1 - \pi_{(i)})}}.
\end{equation}
\end{itemize}
In this paper, we will discuss several variants of HC. To distinguish one from the other, we call the version in (\ref{OHC})  the {\it Orthodox HC (OHC)}.

Fix $0 < \alpha < 1$. To use the HC for a level-$\alpha$ test for the testing problem (\ref{testnull})-(\ref{testalt}),  we must find the critical value $h(p, \alpha)$ defined by
\[
P_{H_0^{(p)}} \{ HC_p^* > h(p, \alpha) \}  = \alpha.
\]
Asymptotically,  it is known that for any fixed $\alpha \in (0,1)$,
\begin{equation} \label{hpalpha}
h(p, \alpha)  = \sqrt{2 \log\log(p)} (1 + o(1)),
\end{equation}
so an approximation of $h(p, \alpha)$ is $\sqrt{2 \log \log(p)}$.  We say $\alpha = \alpha_p$ {\it tends to $0$  slowly enough} if
\[
h(p, \alpha_p)  \sim \sqrt{2 \log \log(p)}.
\]
Consider the HC-test where we reject $H_0^{(p)}$ if and only if
\[
HC_p^* \geq h(p, \alpha_p).
\]
Recall that $\rho^*(\vartheta)$ is the standard phase function defined in \eqref{phasefunction}.
The following theorem is proved by  Donoho and Jin (2004) \cite{DJ04}.
\begin{thm} \label{thm:OHC}
Fix $(\vartheta, r)$ in the phase space such that $r > \rho^*(\vartheta)$.  Suppose as $p \goto \infty$, the level $\alpha_p$ of the HC-test tend to $0$ slowly enough, then the power of the HC-test tends to $1$.
\end{thm}

Combining this with Theorem \ref{thm:Detect} (not requiring $\Omega = I_p$),  for any fixed $(\vartheta, r)$ in Region of Detectable (see (\ref{Detectable})),   OHC yields an asymptotically full power test when $\Omega = I_p$.  Therefore, OHC    achieves the optimal phase diagram for detection.

\begin{remark}
The approximation of  $h(p,\alpha)$  in (\ref{hpalpha}) is largely asymptotic, so it is not very accurate for finite    $p$; our recommendation is to approximate $h(p, \alpha)$ by simulations.
\end{remark}

\begin{remark}
For small $i$, $HC_{p,i}$ is slightly heavy-tailed, so OHC is also heavy-tailed.
To alleviate the problem,  we recommend  the following modified version:
\begin{equation} \label{HC+}
HC_p^{+} =  \max_{\{  1 \leq i \leq \alpha_0 p:  \;     \pi_{(i)}  > 1/p \}}  HC_{p, i}.
\end{equation}
Similarly, we can define $h^+(p, \alpha)$ as the critical value such that
$P_{H_0^{(p)}}(HC_p^+ \geq h(p, \alpha)) = \alpha$.
See Donoho and Jin (2004) \cite{DJ04} for details.
\end{remark}

\begin{remark}
One advantage of OHC-test is that we only need $P$-values to use it,
without any knowledge of the parameters $(\eps_p, \tau_p)$, so the test is not tied to the specific model in (\ref{testnull})-(\ref{testalt}).  On the other hand, in the idealized case where $(\eps_p, \tau_p)$ are known, the optimal test is the
Neyman-Pearson Likelihood Ratio  Test (LRT).  There is an interesting phase transition associated with the limiting distribution of LRT. Write the log-likelihood ratio associated with (\ref{testnull})-(\ref{testalt}) as
\[
LR_p(\eps_p, \tau_p) = LR_p(\eps_p, \tau_p; Y) =   \sum_{i = 1}^p \log((1 - \eps_p) + \eps_p e^{\tau_p Y_i - \tau_p^2/2}).
\]
With the calibrations in (\ref{Defineeps})-(\ref{Definetau}),
$LR_p$ can have non-degenerate limits only when $(\vartheta, r)$ fall {\it exactly} onto the phase boundary $r=\rho^*(\vartheta)$; still, this alone is inadequate, and we must modify the calibrations slightly. In detail,
we let
\begin{equation} \label{Newcalibrations}
r = \rho^*(\vartheta),   \;\;\;  \tau_p = \sqrt{2 r \log(p)},  \;\;\;
\eps_p  =
\left\{
\begin{array}{ll}
p^{-\vartheta},   &\mbox{if $1/2 < \vartheta \leq 3/4$},  \\
\tau_p^{2 \sqrt{r}} p^{-\vartheta},   &\mbox{if $3/4 < \vartheta < 1$}.
\end{array}
\right.
\end{equation}
As $p \goto \infty$, if (\ref{Newcalibrations}) holds, then $LR_p$ has weak limits as follows \cite{JinThesis}:
\[
LR_p  \longrightarrow  \left\{
\begin{array}{ll}
N(\mp \frac{1}{2}, 1), &\mbox{if $1/2 < \vartheta < 3/4$},   \\
N(\mp \frac{1}{4}, \frac{1}{2}), &\mbox{if $\vartheta = 3/4$}, \\
\nu_{\mp}^{(\vartheta)},   &\mbox{if $3/4 < \vartheta < 1$},
\end{array}
\right.
\]
under $H_0^{(p)}$ and $H_1^{(p)}$, respectively.
Here, $\nu_{\mp}^{(\vartheta)}$ are the distributions with the characteristic functions $\psi_{\mp}^{(\vartheta)}$  given by
$\psi_{-}^{(\vartheta)}(t) =
\frac{1}{\sqrt{2 \pi}} \int_{-\infty}^{\infty} [e^{ \sqrt{-1} t \log(1 + e^z)} - 1 - \sqrt{-1} t e^{z}]  e^{-\frac{z}{\vartheta}(1+\sqrt{1-\vartheta})} dz$ and
$\psi_+^{(\vartheta)}(t) =
\frac{1}{\sqrt{2 \pi}} \int_{-\infty}^{\infty} [e^{ \sqrt{-1} t \log(1 + e^z)} - 1]  e^{-\frac{z}{\vartheta}(1-\vartheta+\sqrt{1-\vartheta})} dz$, respectively.
\end{remark}

\subsection{Optimal signal recovery by Hard Thresholding (white noise)}
\label{subsec:recover}
We now consider  the problem of signal recovery in the case of $\Omega = I_p$.  In this simple setting,
\[
Y_i \stackrel{iid}{\sim}  (1 - \eps_p) N(0,1) + \eps_p N(\tau_p,1),
\]
and a conventional approach to estimating  $\beta$ is to use Hard Thresholding (HT).
Fix a threshold $t > 0$.  The HT estimator $\hat{\beta}_t^{HT}$
is given by
\begin{equation*}
\hat{\beta}_{t,i}^{HT}=
   \begin{cases}
 Y_i,  &\qquad \mbox{if  $|Y_i|\geq t$},    \\
 0,     &\qquad \mbox{otherwise}.
\end{cases}
\end{equation*}
It is convenient to choose thresholds having the form of
\[
t_q(p)  =\sqrt{2q\:\log\:p}, \qquad  \mbox{where $0<q<1$ is a fixed parameter}.
\]
Ideally, when $(\vartheta, r)$ are given, we choose $q$ as follows:
\begin{equation} \label{Idealq}
q^{ideal} =  \left\{
\begin{array}{ll}
(\vartheta + r)^2/(4r),  &\qquad  r > \vartheta, \\
\vartheta,  &\qquad 0 < r < \vartheta,
\end{array}
\right.
\end{equation}
and let  $t^{ideal}_p = \sqrt{2 q^{ideal} \log(p)}$.
By direct calculations and Mills' ratio \cite{Wassermanbook},
\[
h_p(\hat{\beta}_{t^{ideal}_p}^{HT}, \beta)    \left\{
\begin{array}{ll}
=  L_p p^{1 - (\vartheta + r)^2/(4r)},  &\qquad \mbox{if $r  > \vartheta$},  \\
\sim  p^{1 - \vartheta},  &\qquad \mbox{if $0 < r < \vartheta$}.
\end{array}
\right.
\]
Combining this with Theorem \ref{thm:recover} (which is for more general $\Omega$), we conclude that when $\Omega = I_p$,   HT achieves the minimax Hamming distance, up to a multi-$\log(p)$ factor;  so it achieves the optimal phase diagram for recovery given in Section \ref{subsec:phase}.

\begin{remark}
The ideal choice of $q$ in (\ref{Idealq})  depends on $(\vartheta, r)$ and  it is hard  to set them in a data driven fashion.  On the other hand, a convenient choice is $q = 1$,   corresponding to the universal threshold $t_p^* = \sqrt{2 \log(p)}$ \cite{Wassermanbook}.  Note that when $r > \rho_{exact}^*(\vartheta; I_p)$,   $\hamm_p(\hat{\beta}_{t_p^*}^{HT}, \beta) \leq C (\log(p))^{-1/2}$ (and so exact recovery is achieved).
\end{remark}

\begin{remark}
When a method yields exact recovery with overwhelming probabilities, we say it has the {\it Oracle Property}, a well-known notion in the  variable selection literature \cite{FanLi}.
 In such a framework, we are using $P(\sgn(\hat{\beta}) \neq \sgn(\beta))$ as the measure of loss.
Seemingly, such a measure is only appropriate for Rare/Strong signals. When signals are Rare and Weak, exact recovery is usually impossible, and the Hamming distance is a more appropriate measure of loss.
\end{remark}

\subsection{Applications} \label{subsec:applications}

HC (and its variants) has found applications in  GWAS  and DNA Copy Number Variation (CNV), where the genetic effects are believed to be rare and weak.
Parkhomenko {\it et al}.  \cite{GWAS1} used HC to detect modest genetic effects
in a genome-wide study of rheumatoid arthritis. Sabatti {\it et al}.  \cite{Sabatti}  used HC
to quantify the strength of the overall genetic signals for each of the nine traits they were interested in.  De la Cruz {\it et al}. \cite{delaCruz}  used HC to test
whether there are associated markers in a given set of markers,
with applications to  Crohn's disease.
Jeng et al. \cite{Jeng1, Jeng2}
proposed a variant of HC called
{\it Proportion Adaptive Segment Selection (PASS)},
which can be viewed as a two-way screening process (across different SNPs and across
different subjects), simultaneously
dealing with the rare genetic effects and rare genomic variation.
See also \cite{GWAS2, Martin, Lin, RW, Wu}.

HC has also found applications in several modern experiments in Cosmology and Astronomy---another
important source of rare and weak signals.  Jin {\it et al}    \cite{Starck} and Cayon {\it et al}   \cite{Cayon1} (see also \cite{Cruz, Cayon2})  applied
HC to standardized wavelet coefficients of Wilkinson Microwave Anisotropy Probe (WMAP),
addressing the problem nonGaussianity detection in the Cosmic Microwave Background (CMB).
Compared to the widely used kurtosis-based non-Gaussianity detector, HC
showed superior power and sensitivity, and  pointed
to the {\it  cold spot}  centered at galactic coordinate  (longitude, latitude) = $(207.8^{\circ}, -56.3^{\circ})$ (see  \cite{Vielva} for more discussions on HC and the cold spot).
Pires {\it et al}  \cite{Pires}  applied many nonGaussianity detectors to gravitational weak lensing data  and showed that
HC is competitive,  being more specifically focused on excess of observations in the tails
of the distribution.
Bennett {\it et al}.  \cite{Bennett} applied the HC ideas to the
problem of Gravitational Wave detection  for a monochromatic periodic source in a binary system. They use a modified form of HC which offers a noticeable increase in the detection power, and yet is robust.

HC has also been applied to disease surveillance for early detection of disease outbreak (Neill and Lingwall \cite{Neill,Neill0})  and to local anomaly detection in a graph (Saligrama and Zhao \cite{anomaly}), where it is found to be have competitive powers.

\subsection{Connections and extensions}
\label{subsec:extensions}
In model (\ref{testnull})-(\ref{testalt}),  the noise entries are iid samples from a distribution $F$ which is known to be $N(0,1)$. Delaigle and Hall \cite{Delaigle} and Delaigle {\it et al}. \cite{Robustness}
address the more realistic setting where $F$ is unknown and is probably nonGaussian.
They consider a two group model (a control group and a case group) and compute $P$-values for each
individual features using bootstrapped Student's $t$-scores.  The problem is also addressed by
 \cite{Greenshtein, ShaoQM}, using very different approaches.

The testing problem (\ref{testnull})-(\ref{testalt}) is a special case of the problem of testing
$H_0^{(p)}$ of  $X_i \stackrel{iid}{\sim}  F$ versus
$H_1^{(p)}$ of  $X_i \stackrel{iid}{\sim} (1 - \eps) F  + \eps  G$,
where $\eps \in (0,1)$ is small,  $F$ and $G$ are two different distributions, and  $(\eps, F, G)$ may depend on $p$.   Cai {\it et al}. \cite{CJJ}  considers the case where $F = N(0,1)$ and $G = N(\tau, \sigma^2)$; see also  \cite{Bogdan1, Bogdan2}.  Park and Ghosh  \cite{Ghosh}  gave a nice review on large-scale multiple testing with a detailed discussion on HC. Cai and Wu \cite{CaiWu} considers the extension where $F = N(0,1)$ and $G$ is a Gaussian location mixture with a general mixing distribution,  and Arias-Castro and Wang \cite{WangMeng} investigate the case where
$F$ is {\it unknown} but symmetric.  In a closely related setting, Laurent {\it et al}.   \cite{Laurent}  considers the problem of testing whether the samples $X_i$ are iid samples from a
 single normal, or a mixture of two normals with different (but unknown) means. Addario-Berry {\it et al}.  \cite{Lugosi}  and Arias-Castro {\it et al}.   \cite{Erycluster}  consider a setting similar to (\ref{testnull})-(\ref{testalt})  but where the signals are structured,  forming clusters in (unknown) geometric shapes; the work is closely related to that in \cite[Section 6]{HJ09}.

Gayraud and Ingster \cite{Ingstervariable} show that HC statistic continues to be successful in detecting very sparse mixtures in a functional setting.
Haupt {\it et al}.  \cite{Haupt1, Haupt2}  consider the setting
where adaptive sample scheme is available so that  we can  do inference and collect
data in an alternating order.

HC can also be viewed as a measure for goodness-of-fit.
Jager and Wellner  \cite{Wellner2007} introduced a new family of goodness-of-fit tests based on the $\phi$-divergence, including HC as a special case.
Wellner and Koltchinskii \cite{WellnerKoltchinskii}  further investigated the Berk-Jones statistic, which is closely related to HC, and derive the limiting distribution of the Berk-Jones statistic.
In Jager and Wellner (2004) \cite{Wellner2004},  they further investigated the limiting distribution of a class of weighted Kolmogorov statistics, including HC as a especial case.
The pontogram of Kendall and Kendall   \cite{Kendall}
is an instance of HC, applied to a special set of $P$-values.

HC is closely related to the False Discovery Rate (FDR) controlling procedure  by  Benjamini and Hochberg (BH)  \cite{BH95}, but it is different in important ways.
The BH procedure targets on Rare/Strong signals,
and the main goal to select the few strong signals embedded in a long list of noise,
without making too many false selections.
HC targets the more delicate Rare/Weak regime, where the signals and noise are hard to distinguishable. In such settings, while
the BH procedure still controls the FDR, it yields very few discoveries. In this case,
a more reasonable goal is to test whether any signals exist without  demanding that we properly identify them all; this is what HC is specifically
designed for.

HC is also intimately connected to the problem of constructing confidence bands for the
{\it False Discovery Proportion} (FDP).  See  \cite{CJL, Una, Ge}. Motivated by a  study of  Kuiper Belt Objects (KBO) (e.g., \cite{Rice}),  Cai {\it et al}.
\cite{CJL} develop HC into an estimator for the proportion of non-null effects, a problem that has
attracted substantial attention in the area of large-scale multiple testing in the past decade.  The literature along this line connects to \cite{BH95} on
controlling FDR, as well as  Efron \cite{EfronNull} on
controlling the local FDR in gene  microarray studies.

\section{Detecting and recovering signals in colored noise}
\label{sec:color}
In this section, we extend the discussions in Section \ref{sec:white} to the case of $\Omega \neq I_p$.
For optimal procedures in the current case, the key is to exploit the graphic structure
of $\Omega$.  We propose Innovated Higher Criticism (IHC) for signal detection and Graphlet Screening (GS)
for signal recovery. IHC and GS  can be viewed as $\Omega$-aware Higher Criticism and $\Omega$-aware Hard Thresholding, respectively.

\subsection{Innovated Higher Criticism and its optimality in signal detection}  \label{subsec:IHC}
We revisit the testing problem (\ref{testnull})-(\ref{testalt}),  where we recall that
\begin{equation} \label{Stein2}
Y = \beta + z, \qquad z \sim N(0, \Sigma).
\end{equation}
Recall that HC is a method of combining  $P$-values.  We are interested in adapting HC for general sparse precision matrix $\Omega$, and there are three perceivable ways  of  combining the $P$-values.

In the first one, we obtain individual $P$-values marginally in a  brute-force fashion:
\[
\pi_i = P(|N(0,1)| \geq |Y_i|  /   (\Sigma(i,i))^{1/2} ).
\]
We call the HC applied to these $P$-values the {\it Brute-force HC (BHC)}. BHC neglects the correlation structure, so we expect some room  for improvement.

For an alternative, denoting the unique square root of $\Omega$ by $\Omega^{1/2}$,
 it is tempting to use the {\it Whitened Transformation}
$Y \mapsto  \Omega^{1/2} Y \sim N(\Omega^{1/2} \beta, I_p)$,
so that the noise is whitened.
We then obtain individual $P$-values by
\[
\pi_i = P(|N(0,1)| \geq   |(\Omega^{1/2} Y)_i|), \qquad 1 \leq i \leq p.
\]
We call the resultant HC the {\it Whitened HC}.

Our proposal is {\it Innovated HC (IHC)}.  Underlying IHC is the idea to find a transformation $Y \mapsto M Y = M\beta + M z$ ($M = M_{p,p}$, may depend on $\Omega$)  so that
\begin{itemize}
\item {\it Preserving sparsity}: most entries of the vector $M \beta$ are zero; this is important since the strength of HC lies in detecting very sparse signals.
\item {\it Simultaneously maximizing SNR}:  to maximize the Signal-to-Noise Ratio for all $i$ at which $\beta_i \neq 0$,  defined by $(M \beta)_i / \sqrt{(MM')(i,i)}$ (since $M Y \sim N(M \beta, MM')$).
\end{itemize}
The best choice  turns out to be $M = \Omega$, associated with which is   Innovated Transformation (IT)
\[
Y  \mapsto   \Omega Y  \sim N(\Omega \beta, \Omega).
\]
This is related to the notion of {\it innovation}  in time series literature and so the name of IT.
See \cite[Section 1.2]{FJY} for detailed discussion on  why $M = \Omega$ is the best choice.

Now,  first, IT preserves the sparsity of $\beta$. Second, for most $i$ at which $\beta_i \neq 0$,  among three choices of $M$, $M = I_p$ (corresponding to model \eqref{Stein2}), $M = \Omega^{1/2}$, and $M = \Omega$, the SNR are
\[
(\Sigma(i,i))^{-1/2} \beta_i,  \qquad ((\Omega^{1/2})(i,i))  \beta_i,  \qquad  \beta_i
\]
correspondingly, where in the last term, we have used the assumption of $\Omega(i,i) = 1$.  See \cite{HJ09} for the insight underlying these results and proofs, where the key is to combine the sparsity of $\Omega$ and the ARW model of $\beta$.  Note that by basic algebra,
\[
\Sigma(i,i)^{-1/2}  \;  \leq \;  (\Omega^{1/2})(i,i)   \;  \leq  \; 1,
\]
so IT has the largest SNR, simultaneously at all $i$ such that $\beta_i \neq 0$.

It is particularly interesting that, while WT yields uncorrelated noise, it does not yield the largest possible SNR, so WT is not the best choice.   For the current setting where signals are Rare and Weak and $\Omega$ is sparse,  larger SNR out-weights sparse correlations among the noise, so we prefer IHC to   WHC. Similarly, we prefer WHC to BHC.

\begin{example} \label{exm:block}
Suppose $\Omega$ is block-wise diagonal and satisfies
$\Omega(i,j) = 1\{i = j\} + h_0 \cdot 1\{|i-j| = 1,  \mbox{$\max\{i,j\}$ is even}\}$,  $-1 < h_0 < 1$,  $1 \leq i, j \leq p$.
For all $1 \leq i \leq p$,   $(\Sigma(i,i))^{-1/2} = \sqrt{1 - h_0^2}$ and $(\Omega^{1/2})(i,i) =  \frac{1}{2} [\sqrt{1 + h_0} + \sqrt{1 - h_0}]$, and so $(\Sigma(i,i))^{-1/2} \leq  (\Omega^{1/2})(i,i)  \leq 1$;
IT yields larger SNR than that of WT, and WT yields larger SNR than  that of model (\ref{Stein2}).
\end{example}

\begin{remark}
At the heart of IHC is entry-wise thresholding applied to the vector $\Omega Y$. This is
equivalent to the Univariate Screening (US) \cite{FanLv, GJW}. In detail, we can rewrite model (\ref{Stein2}) as a regression model $W \sim N(X \beta, I_p)$,  with $W = \Omega^{1/2} Y$ and $X = \Omega^{1/2}$.
US  thresholds the vector $X' W$ entry-wise;  note $X' W = \Omega Y$.
\end{remark}

Similarly,  IHC consists of three simple steps (the last two are the same as in OHC).
\begin{itemize}
\item Obtain two-sided $P$-values by $\pi_i = P(|N(0,1)| \geq |(\Omega Y)_i|)$,  $1 \leq i \leq p$.
\item Sort $P$-values: $\pi_{(1)} < \pi_{(2)} < \ldots < \pi_{(p)}$.
\item Innovated Higher Criticism statistic is then
$IHC_p^*   =  \max_{\{ 1 \leq i \leq p/2 \}}  IHC_{p, i}$,  where $IHC_{p, i}  =  \sqrt{p}  [(i/p)  - \pi_{(i)}]/\sqrt{\pi_{(i)} (1 - \pi_{(i)})}$.
\end{itemize}

Consider the IHC test where we reject
$H_0^{(p)}$ if and only if $IHC_p^* \geq d_p^*(\Omega)  h(p, \alpha)$, where $d_p^*(\Omega)$ is as in (\ref{degreeOmega}) and
$h(p, \alpha)$  is as in Section \ref{subsec:HCbasics}.  It is seen that
\[
P_{H_0^{(p)}}(\mbox{reject $H_0^{(p)}$}) \leq \alpha.
\]
Moreover, we have the following theorem, which extends Theorem \ref{thm:OHC} from the case of $\Omega = I_p$ to the case of more general $\Omega$.
\begin{thm} \label{thm:IHC}
Fix $(\vartheta, r)$ in the phase space.  If $r > \rho^*(\vartheta)$ and as $p \goto \infty$,  $\alpha = \alpha_p$ tends to
$0$ slowly enough, then the power of the IHC-test tends to $1$.  If $r < \rho^*(\vartheta)$, then for any tests, the sum of Type I and Type II testing errors tends to $1$ as $p \goto \infty$.
\end{thm}
Combining this with Theorem \ref{thm:Detect}, for any $(\vartheta, r)$ in the Region of Detectable, IHC provides an asymptotically full power test, so it achieves the optimal phase diagram for detection given in Section \ref{subsec:phase}. 

The proof of Theorem \ref{thm:IHC} has two new ingredients, additional to that of Theorem \ref{thm:OHC}. The first ingredient is Lemma \ref{colornumber} in Section \ref{subsec:detectLB}.
The second one is the similarity between $\beta$ and $\Omega \beta$.
Note that the most interesting region for signal detection is when $1/2 < \vartheta < 1$ \cite{DJ04}. For $\vartheta$ in this range, $\eps_p \ll 1/\sqrt{p}$, so
 $\beta$ has about $p \eps_p$
nonzeros, each equals $\tau_p$, and $\Omega \beta$ has $\lesssim  d_p^*(\Omega)\cdot p \eps_p$
nonzeros, where about $p \eps_p$ of them equals to $\tau_p$, and all others do not exceed $\tau_p$ in magnitude.
Since $d_p^*(\Omega)$ does not exceed a multi-$\log(p)$ term,
we do not expect any difference between the detection boundary of $\Omega$ and that of $I_p$; both are  $r = \rho^*(\vartheta)$.
\subsection{Graphlet Screening and its optimality in variable selection}  \label{subsec:GS}
For signal recovery, we rewrite model (\ref{Stein1}) as a linear regression model
\begin{equation} \label{Stein3}
W  \sim  N(X \beta, I_p), \qquad W  \equiv \Omega^{1/2} Y,  \qquad X \equiv \Omega^{1/2}.
\end{equation}
Our proposal is to use Graphlet Screening (GS), an idea developed recently  by
\cite{JZZ, Ke}.
At the heart of GS is {\it graph-guided multivariate screening}, which intends to overcome the shortcomings of the well-know method Univariate Screening (US) (also called marginal regression \cite{GJW}, Sure Screening \cite{FanLv}), without much increase in computational complexity.
Write
\[
X = [x_1, x_2, \ldots, x_p].
\]
To use US, we project $W$ to the columns of $X$, one at a time, and then apply Hard Thresholding:
\[
\hat{\beta}_j =  (x_j, W)  \cdot  1\{ |(x_j, W)| \geq t\}, \qquad 1 \leq j \leq p, \qquad \mbox{$t > 0$: a threshold}.
\]
The major challenge US faces is {\it signal cancellation} \cite{RW}, meaning that due to the correlations among the design variables,  the SNR of $(x_j, W)$ (which is $E[(x_j, W)] = \sum_{\ell} (x_j, x_{\ell}) \beta_{\ell}$  in the current setting) could be substantially smaller in magnitude than that in the case of orthogonal design.

To alleviate `signal cancellation', an alternative is to use Exhaustive Multivariate Screening (EMS).  Fixing a small integer $m_0$, for any $1 \leq m \leq m_0$ and a subset $\{i_1, i_2, \ldots, i_m\}$,  $1 \leq i_1 < i_2 < \ldots < i_m \leq p$,
we project $W$ to $\{x_{i_1}, x_{i_2},  \ldots, x_{i_m}\}$, and conduct a screening using a $\chi^2$-test with $df = m$.
Unfortunately, EMS is neither computationally feasible nor efficient:  it screens $\sum_{m = 1}^{m_0} \binom{p}{m}$ different subsets of variables, which is hard to handle computationally; also, when we include too many subsets for screening,  we need signals stronger than necessary for successful screening.

GS recognizes that in EMS,  many subsets of variables can be safely skipped for screening, and the key innovation is to use a sparse graph to guide the screening. In detail, let $\cg = (V, E)$ be the graph where $V = \{1, 2, \ldots, p\}$ and there is an edge between nodes $i$ and $j$ if and only if  $\Omega(i,j) \neq 0$.
Let
\[
{\cal A}(m_0) = \{\mbox{all connected subgraphs of $\cg$ with size $\leq m_0$} \}.
\]
Compared to EMS, the difference is that GS only applies $\chi^2$-screening to those subsets in
${\cal A}(m_0)$.  Note that when $\Omega = I_p$,  GS reduces to  Hard Thresholding, so it can be viewed as an $\Omega$-aware Hard Thresholding.

By a well-known result in graph theory \cite{Frieze},
\begin{equation} \label{size}
|{\cal A}(m_0)| \leq C p (e d_p^*(\cG))^{m_0},
\end{equation}
where $d_p^*(\cG)$ is the maximum degree of $\cG$. By the definition and (\ref{degreeOmega}),  $d_p^*(\cG) = d_p^*(\Omega)$, and does not exceed a multi-$\log(p)$ term.
As a result, GS has a much smaller computational cost  than that of EMS (in fact, it is only larger than that of US by a multi-$\log(p)$ factor for fixed $m_0$), and also requires much weaker signals than  EMS does for successful screening.

\begin{remark}
GS is a flexible idea and can be adapted to many different settings, where the implementation
may vary from occurrence to occurrence. It is a screening method and it has been applied to variable selection \cite{JZZ, Ke}, which includes model (\ref{Stein3}) as a special case.  It can also been viewed as a way to evaluate the combined significance of (a small number of) features, so it can be used for feature ranking; see Section \ref{sec:future} for more discussions.
\end{remark}

We now describe how to apply GS to model (\ref{Stein3}) for signal recovery.
List all elements in ${\cal A}(m_0)$ in the order of sizes, with ties breaking
lexicographically,
\[
\call_1, \call_2, \ldots, \call_N, \qquad N \equiv |{\cal A}(m_0)|.
\]
Our proposal is a two-step procedure, containing a Screen step and Clean step.
Fix positive tuning parameters $(u, v, q)$.
In the Screen step,  initialize with ${\cal S}_0 = \emptyset$.  For $i = 1, 2, \ldots, N$, letting ${\cal S}_{i-1}$ be the set of all retained indices up to  stage $i-1$,  we update ${\cal S}_{i-1}$ by
\[
{\cal S}_i  =  \left\{
\begin{array}{ll}
{\cal S}_{i-1} \cup \call_i, &   \mbox{if $\| P^{\call_i} Y \|^2  - \|P^{\call_i \cap {\cal S}_{i-1}} Y \|^2 \geq 2 q \log(p)$}, \\
{\cal S}_{i-1}, &   \mbox{otherwise}
\end{array}
\right.
\]
where for any $\call \subset \{1, 2, \ldots, p\}$, $P^{\call}$ is the projection matrix from $R^n$ to $\{x_j: j \in \call\}$.
The set of all retained nodes in the Screen step is then ${\cal S}_N$.

In the Clean step,   when $j \notin  {\cal S}_N$, we set $\hat{\beta}_j^{gs} = 0$.
When $j  \in {\cal S}_N$, let $\cG_{{\cal S}_N}$ be the subgraph of $\cG$ formed by restricting all nodes to ${\cal S}_N$.  We decompose
\[
\cG_{{\cal S}_N} =  \cG_{{\cal S}_N, 1} \cup \cG_{{\cal S}_N, 2} \cup \ldots \cup \cG_{{\cal S}_N, L},
\]
and estimate $\{\beta_j: j \in \cG_{{\cal S}_N, \ell}\}$, $1 \leq \ell \leq L$,  by minimizing
\[
\| P^{\cG_{{\cal S}_N, \ell}} (Y -  \sum_{j \in \cG_{{\cal S}_N, \ell}} \beta_j x_j) \|^2 + u^2 \| \beta \|_0,
\]
subject to the constraint that either  $\beta_j = 0$  or  $|\beta_j| \geq v$.  Putting these together gives the final estimate, denoted by $\hat{\beta} = \hat{\beta}^{gs}(m_0, u, v, q)$.

\begin{thm}  \label{thm:entrywise}
Fix $(m_0, \vartheta, r)$ such that
$1 <  r /  \vartheta <  3 + 2 \sqrt{2} \approx 5.828$ and $m_0 \geq (r-\vartheta)^2/(4\vartheta r)$.
Suppose  (\ref{diagOmega})-(\ref{degreeOmega}) hold, the spectral norm of $\Omega^{-1}$ is bounded by $C$, and $\max_{1\leq i\leq p}\sum_{j=1}^p|\Omega(i,j)|^{\gamma}\leq C$, for some constants $\gamma\in (0,1)$ and $C>0$.
Also, suppose
$|\Omega(i,j)|  \leq    4 \sqrt{2} - 5 \approx 0.6569$  for all $1 \leq i, j \leq p$, $i \neq j$.  If we set the tuning parameters $(u,v, q)$ in GS by $u=\sqrt{2\vartheta\log(p)}$, $v=\sqrt{2r\log(p)}$ and $q$ an appropriately small constant,  then as $p \goto \infty$,
$h_p(\hat{\beta}^{gs}, \beta) \leq  L_p    \hamm_p^*(\vartheta, r; \Omega) =  L_p p^{1 -(\vartheta + r)^2/(4r)} + o(1)$.
\end{thm}
Recall that $\hamm_p^*(\vartheta, r; \Omega)$ is the minimax Hamming distance as in (\ref{eqhamm}).
This says that  for all $\Omega$ considered in Theorem \ref{thm:entrywise}, GS achieves the optimal phase diagram for  recovery; see Section \ref{subsec:phase}.

Theorem \ref{thm:entrywise} was proved in Jin {\it et al} (2012)  \cite[Section 2.6]{JZZ}, as a special example.
The conditions on $r/\vartheta$ and off-diagonals of $\Omega$ are not necessary for GS to achieve the minimax Hamming distance. In fact,  $h_p(\hat{\beta}^{gs}, \beta) \leq  L_p    \hamm_p^*(\vartheta, r; \Omega) + o(1)$ holds in much broader settings, but $\hamm_p^*(\vartheta, r; \Omega)$ can not have such a simple expression.
See more discussions in \cite{JZZ} for the asymptotic minimaxity of GS in more general settings.

\begin{remark}
For the tuning parameters, $m_0$ is usually chosen for the computational capacity. The choice of  $q$ is relatively flexible, as long as it falls into certain ranges. The choice of  $u$ is harder, but the best $u$ is a function of $\eps_p$;  in some settings (e.g., \cite{CJL}), we can estimate $\eps_p$ consistently, and we know how to choose the best $u$. For these reasons, we essentially only have one tuning parameter $v$, which is connected to the tuning parameter in the subset selection and that of the lasso.   See \cite{JZZ} for more discussions.
\end{remark}

\subsection{Phase diagrams (colored noise)} \label{subsec:phase2}
Recall that the optimal phase diagram for general $\Omega$ consists of
four subregions separated by three curves $r=\rho^*(\vartheta)$, $r=\vartheta$ and $r=\rho^*_{exact}(\vartheta; \Omega)$;     
 $\rho^*_{exact}(\vartheta;  \Omega)$ may depend on the off-diagonals of $\Omega$ in  a complicated way, but   we always have
\[
\rho_{exact}^*(\vartheta; \Omega) \geq  \rho_{exact}^*(\vartheta; I_p),
\]
since $\Omega=I_p$ is the easiest case for exact recovery.

For $\Omega$ satisfying conditions of Theorem \ref{thm:entrywise} and for $(\vartheta, r)$ such that $1<r/\vartheta<3+2\sqrt{2}$, the minimax Hamming distance for $\Omega$ has the same convergence rate as that for the case of  $\Omega = I_p$. Note that
in the phase space, the curve  $r = \rho_{exact}^*(\vartheta; I_p)$ and the line $r / \vartheta = 3 + 2 \sqrt{2}$ intersect  at the point $(\vartheta, r) = (1/2, (3 + 2 \sqrt{2})/2)$.  Therefore,
\[
\rho^*_{exact}(\vartheta;  \Omega) =  \rho^*_{exact}(\vartheta;  I_p), \qquad \mbox{for all $1/2 < \vartheta < 1$}.
\]
Consequently, the right half of the curve $r = \rho^*_{exact}(\vartheta;  \Omega)$  coincides with the right half of the curve  $r = \rho_{exact}^*(\vartheta;  I_p)$. See \cite{JZZ} for discussion in more general settings.

By Theorems \ref{thm:Detect}-\ref{thm:recover} and Theorems \ref{thm:IHC}-\ref{thm:entrywise},  the optimal phase diagram for detection is achieved by IHC,  and
the optimal phase diagram for recovery is achieved by GS for a wide range of $\Omega$,
 including but are not limited to those  satisfying the  conditions of Theorem \ref{thm:entrywise}.
 See \cite{JZZ} for  details.

\subsection{An example,  and comparisons with $L^0/L^1$-penalization methods}  \label{subsec:example}

In general, it is hard to derive an explicit form for $r =  \rho_{exact}^*(\vartheta; \Omega)$ for the whole range of $\vartheta$. Still, examples for some $\Omega \neq I_p$ would shed light on how this curve depends on the off-diagonal
entries of $\Omega$.

\begin{figure}[tb]
\centering
\includegraphics[width = 5.8  in, height = 2.2 in]{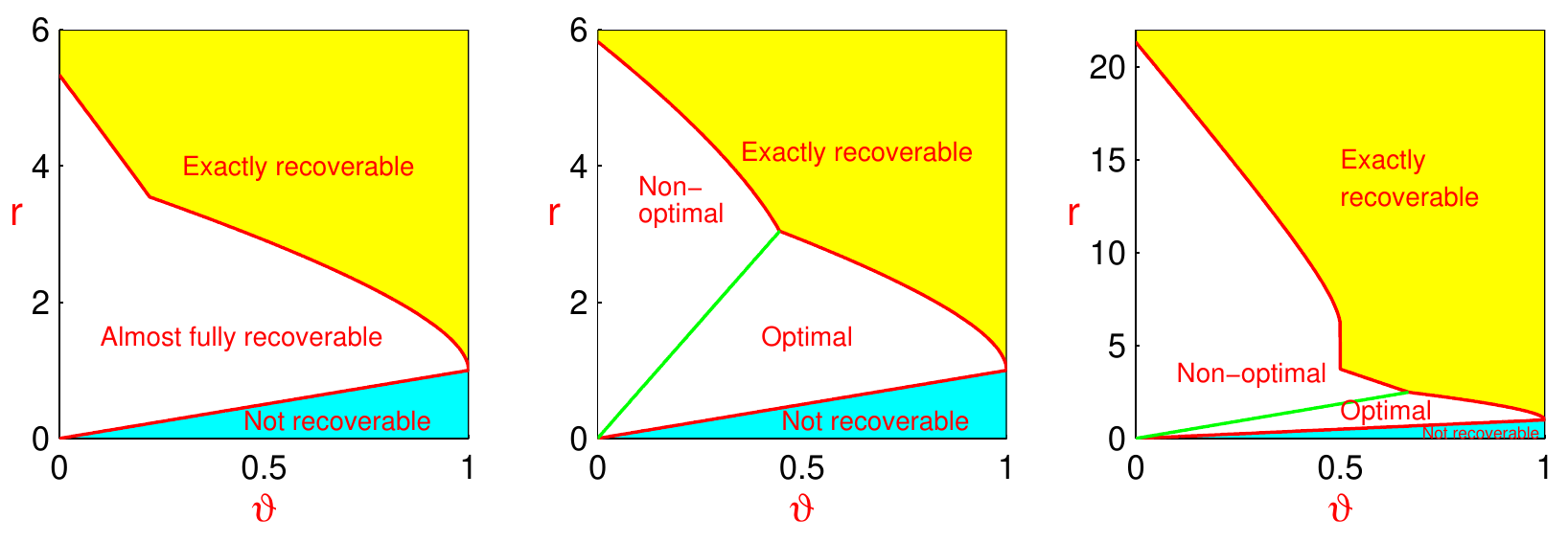}
\caption{Phase diagrams (block-wise diagonal example,  $h_0=0.5$). From left to right:  GS, $L^0$, and $L^1$-penalization method. Note that the first two subregions described in Section \ref{subsec:phase} are combined into Region of Not Recoverable, for convenience.}
\label{fig:bdphase}
\end{figure}

We revisit Example \ref{exm:block} in Section \ref{subsec:IHC}, where $\Omega$ is block-wise diagonal, and each diagonal block is the $2 \times 2$ matrix with $1$ on the diagonals and $h_0$ on the off-diagonals.
It was shown in  \cite{JZZ} that $\hamm^*_p(\vartheta, r,   \Omega)= L_p p^{1 -c(\vartheta,r;  h_0)}$, where
\begin{equation}\label{eq:blockGSrate}
c(\vartheta, r;  h_0) =   \min \biggl\{\frac{(\vartheta +r)^2}{4r},  \;   \vartheta + \frac{(1 - | h_0|)}{2} r,   \;  2 \vartheta +   \frac{\{[(1 - h_0^2) r  - \vartheta]_+\}^2}{4 (1 -  h_0^2) r}  \biggr\}.
\end{equation}
The curve $\rho^*_{exact}(\vartheta; \Omega)$ is then the solution of $c(\vartheta,r;  h_0)=1$, which depends on $h_0$. The top left panel of Figure~\ref{fig:bdphase} displays the phase diagram for $h_0=0.5$.

Somewhat surprisingly, even for very simple $\Omega$ such as the block-wise diagonal example above and even when the tuning parameters are ideally set, subset selection ($L^0$-penalization) and the lasso have non-optimal phase diagrams; in particular, their Region of Exactly Recoverable is smaller than that of GS.
Figure~\ref{fig:bdphase} shows  phase diagrams associate with GS, $L^0/L^1$-penalization methods for the block-wise diagonal example, where the tuning parameters are set ideally to minimize the Hamming distance; see Ji and Jin (2010) \cite{JiJin} and Jin {\it et al.} (2012) \cite{JZZ}. Given the non-optimality of $L^0/L^1$-penalization methods in such a simple $\Omega$, we do not expect that for more general $\Omega$ they could be optimal.

We must note that the optimalities of $L^0/L^1$-penalization methods in the literature are largely limited to settings different from here, where they usually have Rare/Strong signals, and either $\ell^2$-loss or $P(\sgn(\hat{\beta})\neq \sgn(\beta))$
is frequently used as the measure of loss.  However, $\ell^2$-loss is more appropriate for prediction setting, not for variable selection, and $P(\sgn(\hat{\beta})\neq\sgn(\beta))$ is more appropriate for Rare/Strong signals, not for
Rare/Weak signals where it is merely impossible to fully recover the support of $\beta$.
Since $L^0$-penalization method is the target of many penalization methods, including the lasso, SCAD \cite{FanLi},
MC+ \cite{MCP}, we should not expect these penalization methods to be optimal as well.

\begin{remark} Ke {\it et al.} (2012) \cite{Ke}
studied a more complicate setting than that in Theorem \ref{thm:entrywise} or that in \cite{JZZ},  where the Gram matrix is not sparse but is sparsifiable.
They derived the phase diagrams for a case that $\Omega$ is the correlation matrix of a long-memory time series and for a change-point model. The change-point model is a special case of model \eqref{Stein3} where $X$ is an upper triangular matrix of $1$'s (therefore, $X\beta$ is piece-wise constant). For the change-point model, the phase space partitions into only $2$ regions, separating by the curve
$r = \rho_{exact}^{*,cp}(\vartheta)$, where $\rho_{exact}^{*,cp}(\vartheta)  = \max\{4(1 - \vartheta),   (4-10\vartheta)+2\sqrt{[(2-5\vartheta)^2-\vartheta^2]_+} \}$.
\end{remark}

\subsection{Connections to the literature}
Model \eqref{Stein2} is closely related to the linear model $W = N(X\beta, I_n)$, where the Innovated Transformation
reduces that of  $W \mapsto X'W$.
Arias-Castro {\it et al.} (2011) \cite{Castro2} applied HC to $X'W$ for signal detection, which is similar to IHC.
Ingster {\it et al.}  (2009) \cite{IPV} considered a case that $W\sim N(X\beta, \sigma^2 I_n)$ for an unknown $\sigma$ and $X_i$'s are iid samples from $N(0, (1/n)I_p)$. They proposed a modified IT, $W\mapsto \Vert W\Vert^{-1}X'W$, to adapt to the unknown $\sigma$. Mukherjee {\it et al.} (2013) \cite{Lin} considered the binary-response logistic regression. They proposed HC-like statistics for signal detection and exposed interesting dependence of the detection boundary on the design matrix.

Another related setting is that the data are iid samples $Y_1,\cdots, Y_n$ of $N(\beta, \Sigma)$. This  reduces to model \eqref{Stein2} noting that $(1/\sqrt{n})\sum_{i=1}^n Y_i\sim N(\beta, \Sigma)$ is the vector of sufficient statistics of $\beta$. When the data are nonGaussian, Zhong {\it et al.} \cite{Chen} proposed an ``$L_{\gamma}$-thresholding test" which takes BHC as a special case of $\gamma=0$.

GS, as a method to improve US, is different from the Iterative Sure Independence Screening (ISIS) \cite{FanLv,ISIS}. ISIS first applies US to select a small set of variables $M_1$. In the second step, for each $j\notin M_1$, it runs a least-square algorithm  on the model $M_1\cup\{j\}$ and records the coefficient of $j$. These coefficients are then used to rank variables and expand $M_1$ to a set $M_2$. This procedure runs iteratively.
ISIS alleviates `signal cancellation' between variables in $M_1$ and those in $\{1,\cdots, p\}\backslash M_1$, but unlike GS, it does not deal with `signal cancellation' among variables in $\{1,\cdots, p\}\backslash M_1$.

GS is closely related to LARS \cite{LARS} and forward-backward greedy algorithm \cite{FoBo} in utilizing local graphic structure of variables. The Screen step of GS is a step-wise forward algorithm and the Clean step is a backward algorithm.

\section{Stylized applications}  \label{sec:future}
HC and GS are flexible ideas that can be adapted to
a broad set of problems and settings.
In this section, we outline some
potential applications  of HC and GS.

\subsection{Higher Criticism for estimating the bandwidth of a  matrix}
The HC idea, although still in its early stage of development, is seeing increasing interest both in practice and in theory. 
In Section \ref{subsec:applications}-\ref{subsec:extensions}, we have reviewed applications and extensions of HC
in  many different settings.
In this section, we illustrate a new application of HC.

Suppose we have samples $X_i \in R^p$ from a Gaussian distribution
\[
X_i\overset{iid}{\sim} N(0, \Sigma), \qquad 1\leq i\leq n.
\]
The Gaussian assumption is not critical and is only for simplicity.
In many applications, with the Linkage Disequilibrium (LD) matrix  being an iconic example,
$\Sigma$ is unknown but is banded; denote the bandwidth by $b = b(\Sigma)$ so that $b$ is the smallest integer such that
$\Sigma(i,j) = 0$ for all $i,j$ with $|i  - j| \geq b+1$.

We adapt HC to estimate  $b(\Sigma)$. HC can also be adapted
to test whether $b(\Sigma) \leq k_0$ or $b(\Sigma)  >  k_0$ for a given small integer $k_0$; the discussion is similar so we omit it to save space.

Let  the empirical covariance matrix be $S_n = \frac{1}{n}\sum_{i=1}^n X_iX_i'$. For $1 \leq k \leq p-1$, let  $\xi^{(k)}$  and $\hat{\xi}^{(k)}$ be the $(p-k) \times 1$ vectors formed by the  $k$-th (upper) off-diagonal of $\Sigma$ and $S_n$, respectively:
\[
\xi^{(k)}  =  (\Sigma(1,1+k), \ldots, \Sigma(p - k, p))',   \qquad
\hat{\xi}^{(k)}  =  (S_n(1,1+k),   \ldots, S_n(p - k, p))'.
\]
We consider a Rare/Weak setting where each $\xi^{(k)}$ has a small fraction of nonzeros,  and each nonzero is relatively small. By elementary statistics, for any $i,j$ such that $\Sigma(i,j) = 0$, we have that approximately,
$\sqrt{n} S_n(i,j)  \sim  N(0, 1)$.

We propose the following HC estimator for $b(\Sigma)$.
Fix an integer $b_0$ (a relatively small but conservative upper bound for $b(\Sigma)$) and a level $\alpha \in (0,1)$,
\begin{itemize}
\item For $k = 1, \ldots, b_0$,  apply $HC_p^{+}$ in (\ref{HC+})   to $\hat{\xi}^{(k)}$, where the individual
$P$-values associated with the entries of $\hat{\xi}^{(k)}$ are computed by
$P(|N(0,1)| \geq \sqrt{n} \hat{\xi}^{(k)}_i)$, $1 \leq i \leq p-k$.
Denote the resultant HC scores by $HC^{(1)}, \ldots, HC^{(b_0)}$, correspondingly.
\item Estimate $b(\Sigma)$ by
$\hat{b}^{HC} = \hat{b}^{HC}(S_n; n, p, t_n, b_0)  = \max\{1 \leq k \leq b_0: HC^{(k)} \geq h^+(p, \alpha/b_0)\}$.
\end{itemize}
Recall that $h^+(p, \alpha)$ is as in Section \ref{subsec:HCbasics} which can be computed by simulations.

We conducted a small-scale simulation, where $(p, n, b(\Sigma), b_0, \alpha) = (5000,$ $200, 2, 10, 0.05)$. For $k = 1, 2$, and fixed $(\eps, \tau)$, we generate the entries of  $\xi^{(k)}$ randomly from $(1 - \eps) \nu_0 + \eps \nu_{\tau}$.
We then apply the above procedure and repeat the whole simulation processes independently
 for $200$ times, and recorded the error rates (the fraction of simulations where $\hat{b}^{HC} \neq b(\Sigma)$).
We have investigated $6$ different combinations of $(\eps, \tau)$: $(.01, .175)$, $(.01, .2)$, $(.01, .225)$, $(.005, .225)$, $(.005, .25)$ and $(.01, .275)$; and the corresponding error rates of $\hat{b}^{HC}$ are $6.5\%$, $0.5\%$, $0\%$, $8.5\%$, $3\%$ and $2\%$.

\begin{remark}
The choice of $h^+(p, \alpha/N)$ is based on Bonferroni correction which is acceptable for relatively small $N$. For large $N$, we may need to  adjust the threshold, say, with Benjamini and Hochberg's FDR-controlling method \cite{BH95}.
\end{remark}

\subsection{Ranking features by Graphlet Screening} \label{subsec:ranking}
Consider a linear regression of $n$ samples and $p$ features (or variables):
\begin{equation} \label{Stein4}
W = X \beta + z,  \qquad    X = X_{n,p} = [x_1, x_2, \ldots, x_p], \qquad z \sim N(0, I_n).
\end{equation}
Denote the Gram matrix by
\[
G = X'X.
\]
We assume $X$ is normalized so that $G$ has unit diagonals.  We are primarily interested in the case where $G$ is approximately sparse in the sense that each row of $G$ has relatively few large entries.
At the same time, we assume the feature vector $\beta$ is sparse in the sense that it only has a small fraction of nonzeros.
We are interested in ranking the features in a way to have a competitive Receiver Operating Curve Characteristic (ROC) curve.

Conventionally, we rank the features by US: we project $W$ to the columns of $X$, one at a time, and rank their significances according to $|(x_j, W)|$.
The challenge of this approach is, again, `signal cancellation'; see
Section \ref{subsec:GS} for discussions. When `signal cancelation' presents, there  is room for improvement.

We propose to rank the features by GS.  Fixing a threshold $\delta > 0$,
introduce a graph $\cg^{*, \delta} = (V, E)$  where $V = \{1, 2, \ldots, p\}$ and there is an edge
between $i$ and $j$ if and only if $|G(i,j)|\geq \delta$;  since $G$ is approximately sparse,  $\cg^{*,\delta}$ is sparse in that the maximum degree is small, given an appropriate choice of $\delta$.
Fixing $m_0 > 1$,  we similarly define
\[
{\cal A}^{*, \delta}(m_0)={\cal A}^{*,\delta}(m_0, G) = \{ \mbox{all connected subgraphs of $\cg^{*,\delta}$ with size $\leq m_0$}\},
\]
Similarly to \eqref{size},
\[
|{\cal A}^{*,\delta}(m_0)| \leq C p (e d_p^*)^{m_0},   \qquad \mbox{where $d_p^* = d_p^*(\delta, G)$ is the maximal degree of $\cg^{*,\delta}$}.
\]
Our procedure consists of the following steps.
\begin{itemize}
\item For each $\call\in{\cal A}^{*,\delta}(m_0)$, compute a $P$-value by $\pi^{(\call)} = P(\chi^2_{|\call|}(0)>
 \|P^{\call} Y\|^2)$.
\item Define $\pi^{gs}_j = \min_{\call\in{\cal A}^{*,\delta}(m_0)} \{ \pi^{(\call)} \}$, for $1\leq j\leq p$.
\item Rank the significance of feature $j$ according to $\pi^{gs}_j$.
\end{itemize}
Recall that $P^{\call}$ is the projection from $R^n$ to $\{x_j: j \in \call \}$. The procedure is related to the  hierarchical variable selection procedures \cite{Meinshausen} but is
different in significant ways.

\begin{figure}[tb]
\centering
\includegraphics[height = 2 in, width = 2.8 in]{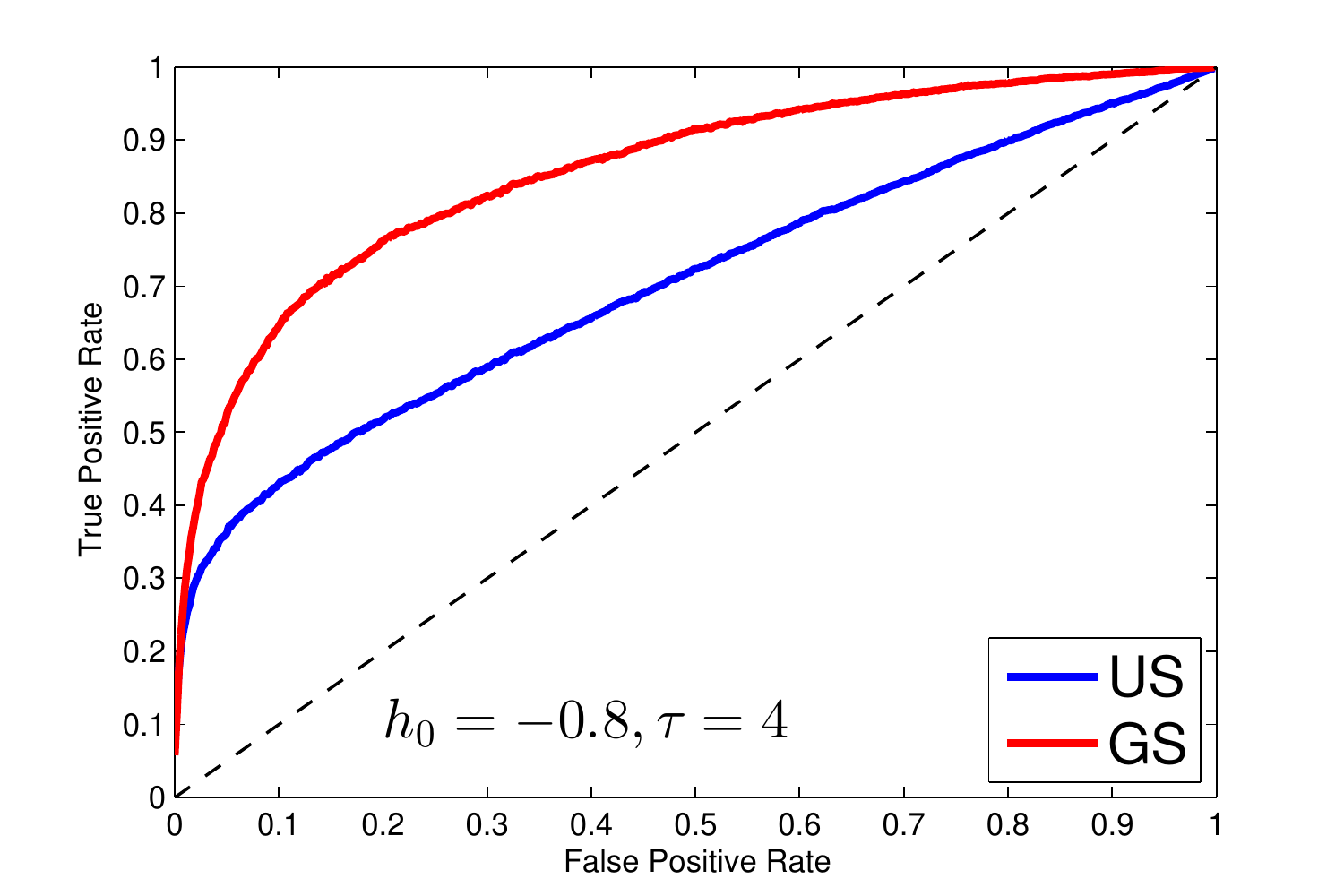}
\includegraphics[height = 2 in, width = 2.8 in]{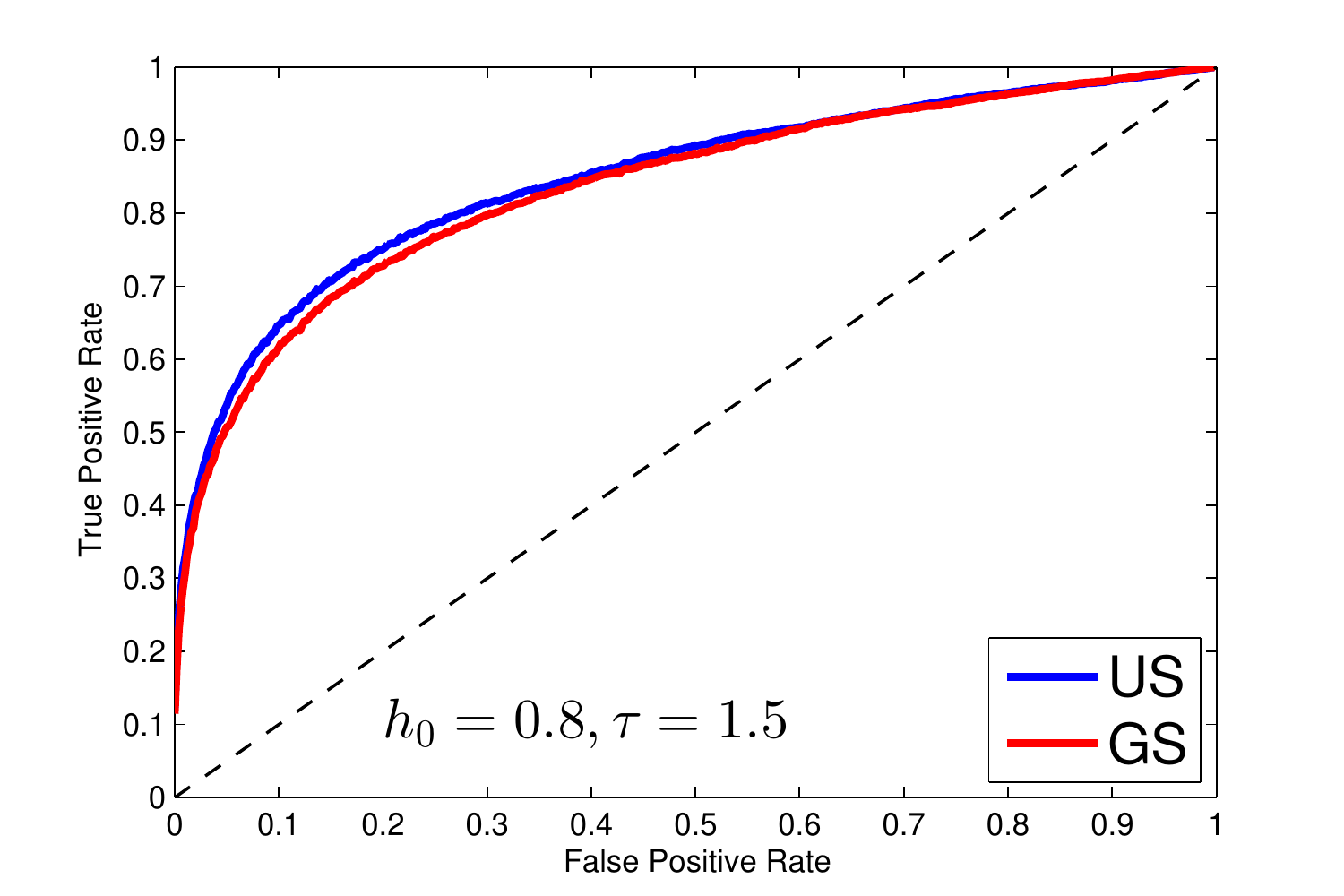}
\caption{ROC curves associated with features ranking by US (blue) and GS (red).
Signal cancellation is severe on the left and GS offers a significant improvement. It is not severe on the right so GS is comparable to US; see Section \ref{subsec:ranking} for details.  }
\label{fig:GSpvals}
\end{figure}

We conducted a small-scale numerical study, where $(n, p,\eps) =(500, 1000, 0.05)$. Let $\Sigma$ be a $p\times p$ blockwise diagonal matrix with size-$2$ blocks; each block has diagonals $1$ and off-diagonals $h_0$. Given $(h_0,\tau)$, we first generate $(\beta_{2j-1},\beta_{2j})\overset{iid}{\sim}(1-\eps)\nu_{(0,0)} + (\eps/2)\nu_{(\tau,\tau)}+(\eps/2)\nu_{(\tau,0)}$, for $j=1,\cdots, p/2$, where $\nu_{a}$ is a point mass at $a$ for any $a\in R^2$. Next, we generate $X_i\overset{iid}{\sim} N(\beta, (1/n)\Sigma)$, for $i=1,\cdots, n$.
We applied both US and GS (taking $m_0=2$) to rank features. Figure~\ref{fig:GSpvals} displays the corresponding ROC curves, obtained from averaging $200$ independent repetitions. We have investigated two cases $(h_0,\tau)=(-0.8,4), (0.8, 1.5)$. In the first case, signal cancellation is severe and GS significantly outperforms US. In the second case,
GS has a similar performance as US.

Feature ranking is of interest in many high dimensional problems,
including but are not limited to (a) large-scale multiple testing, where it is of interest to
develop methods that control the FDR while maximizing the power of multiple tests,
(b) cancer classification where it is desirable to select a small fraction of features for the
trained classification decision \cite{DJ08, DJ09, JinPNAS}, and spectral clustering where it is desirable to perform a dimension reduction before we apply Principle Component Analysis (PCA) \cite{JW}.
As GS provides a better strategy in feature ranking than US, it is potentially useful  in attacking all the problems above.

\section{Feature selection by Higher Criticism for classification}
\label{sec:class}
Among many uses of Higher Criticism, one that is particularly interesting is to set thresholds for
feature selection in the context of  classification.

Consider a two-class classification setting where
$(Y^{(i)}, \ell_i)$,  $1 \leq i \leq n$,   are measurements from two different classes.  Here,
 $Y^{(i)} \in R^p$ are the feature vectors and $\ell_i  \in\{-1, 1\}$ are the class labels.
We assume two classes are equally likely, so that after a standardizing transformation,
\[
Y^{(i)}  \sim  N(\ell_i \cdot  \mu, \Sigma),
\]
with $\mu\in R^p$ the contrast mean vector and $\Sigma$ the $p \times p$ covariance matrix;
such an assumption is only for simplicity in presentations.
Given a fresh feature vector $Y$,  the primary interest is to predict the associated class label $\ell  \in \{-1, 1\}$.

For simplicity,  we assume $\Sigma$ is known and the precision matrix $\Omega = \Sigma^{-1}$ is sparse.
The case $\Sigma$ is unknown (but $\Omega$ is sparse)  is discussed in Fan {\it et al.} (2013) \cite{FJY};
see details therein.

Fisher's linear discriminant analysis (LDA) is a classical approach to classification.
Let $w = (w_1, w_2, \ldots, w_p)' $ be a $p \times 1$  feature weight vector. For a fresh feature vector $Y=(Y_1,\cdots, Y_p)'$, Fisher's LDA takes the form
\[
L(Y)  = \sum_{i =1}^p w_i  Y_i,
\]
and classifies $\ell= \pm 1$ according to  $L(Y) \gtrless 0$.
When $(\Sigma, \mu)$ are known,  it is known that the  optimal weight vector satisfies
$w \propto  \Omega \mu$.

To adapt Fisher's LDA to the current setting,  the key is to estimate $\mu$.
We are primarily interested in the Rare/Weak setting where only a small fraction of the entries of $\mu$ is nonzero and
the nonzero entries are individually small. Define the feature $z$-vector
\begin{equation} \label{Zvector}
Z   =   \frac{1}{\sqrt{n}} \sum_{i = 1}^n  (\ell_i \cdot  Y^{(i)}) \sim N(\sqrt{n} \mu, \Sigma).
\end{equation}
A standard approach to estimating $\mu$ is by some sort of thresholding scheme.
For any $t > 0$,  denote by $\eta_t(z)$ the clipping thresholding function $\eta_t(z) = \sgn(z) 1\{|z| \geq t\}$ \cite{DJ08,FJY}.
Our proposal is to use Innovated Thresholding which thresholds $\Omega Z$ coordinate-wise:
\begin{equation} \label{IT}
\hat{\mu}_{t,i}^{IT} =  \eta_{t}((\Omega Z)_i), \qquad  1 \leq i \leq p.
\end{equation}
One may  also use Brute-force Thresholding which thresholds $Z$ coordinate-wise, or the Whitened Thresholding which thresholds  $\Omega^{1/2} Z$ coordinate-wise. However,
these schemes are inferior to Innovated Thresholding,
for Innovated Transformation yields larger Signal-to-Noise Ratio than  Brute-force Transformation and  Whitened Transformation; see Section \ref{subsec:IHC} for details.
Also,  in (\ref{IT}), we have used the clipping thresholding rule.  Alternatively, one may use hard-thresholding or soft-thresholding, but the difference is usually not significant.  See \cite{DJ08, FJY}
for more discussions.

We now modify Fisher's LDA by
\[
L_{t}^{IT}(Y; \Omega) = (\hat{\mu}_t^{IT})' \Omega Y, \qquad \mbox{where} \qquad  \hat{\mu}_t^{IT} = (\hat{\mu}_{t,1}^{IT}, \hat{\mu}_{t,2}^{IT}, \ldots, \hat{\mu}_{t,p}^{IT})'.
\]
and classify $\ell$ as $\pm 1$ according to $L_t^{IT}(Y; \Omega) \gtrless 0$.
This is related to the modified HC in \cite{Chen}, but  the focus there is on signal detection, not feature selection.

Seemingly, an important issue is how to set the threshold $t$.  We propose
to set the threshold by Higher Criticism Thresholding (HCT), a variant of OHC.
Fix  $\alpha_0 \in (0,1/2]$.
\begin{itemize}
\item Calculate  (two-sided) $P$-values
$\pi_i =  P\{ |N(0,1)| \geq |(\Omega Z)_i|\}$,  $1 \leq i \leq p$.
\item Sort the $P$-values into ascending order:  $\pi_{(1)} < \pi_{(2)} < \ldots  < \pi_{(p)}$.      \item Define the {\it Higher Criticism feature scores} by
\begin{equation}  \label{HCobjclasf}
HC(i;  \pi_{(i)})   =  \sqrt{p}  \frac{i/p   -  \pi_{(i)}}{\sqrt{(i/p)(1 - i/p)}}, \qquad 1 \leq i \leq p.
\end{equation}
Obtain the maximizing index of $HC(i; \pi_{(i)})$:
\[
\hat{i}^{HC} = \margmax_{\{1 \leq i \leq \alpha_0 \cdot p \}} \{HC(i; \pi_{(i)}) \}.
\]
The {\it Higher Criticism threshold (HCT)} for feature selection is then by
\[
\that_p^{HC}   = \that_p^{HC}(Z_1, Z_2, \ldots, Z_p;  \alpha_0)  = |Z|_{\hat{i}^{HC}}.
\]
\end{itemize}
In practice,  we set $\alpha_0 = 0.10$; HCT is  relatively insensitive to  different choices of $\alpha_0$.
In (\ref{HCobjclasf}),  the denominator of the HC objective function  is different from
that of OHC  we used for testing problems (\ref{testnull})-(\ref{testalt}),  although in a similar spirit.
See \cite{DJ09} for explanations.

Once the threshold is decided, the associated Fisher's LDA  is then
\begin{equation} \label{HCTrule}
L_{HC}^{IT}(Y; \Omega) = (\hat{\mu}_{HC}^{IT})' \Omega  Y, \;\;\;     \mbox{where} \;   \hat{\mu}_{HC}^{IT}  = \hat{\mu}_{t}^{IT}   \bigr|_{t = \hat{t}_p^{HC}}.
\end{equation}
The HCT trained classification rule  classifies $\ell = \pm 1$ according to   $L_{HC}^{IT}(Y)  \gtrless  0$.

\begin{remark}
The classification problem is closely connected to the testing problem (\ref{testnull})-(\ref{testalt})  in Sections \ref{sec:white}-\ref{sec:color}. For illustration, assume $\Omega = I_p$ and $\sqrt{n} \mu_j  \stackrel{iid}{\sim} (1 - \eps) \nu_0 + \eps \nu_{\tau}$.
Given a test feature $Y  \sim N(\ell \cdot  \mu, I_p)$,  the classification problem can be viewed as the problem of testing
$H_0^{(p)}$ of $Y \sim N(-\mu, I_p)$ against $H_1^{(p)}$ of $Y  \sim N(\mu, I_p)$.
Despite that this is very similar to that of  (\ref{testnull})-(\ref{testalt}),
there is a major difference. In (\ref{testnull})-(\ref{testalt}), we don't have any information additional to
 the prior distribution on $\mu$,  so all features are equally likely to be useful. In the classification problem, however, 
 the training $z$-vector $Z \sim N(\sqrt{n} \mu, I_p)$ contains additional information about $\mu$; for feature $i$, $1 \leq i \leq p$, the  posterior probability that it is a useful feature  is given by
$P(\mu_i \neq 0 | Z) =  \eps e^{\tau Z_i - \tau^2/2} / [(1 - \eps) + \eps e^{\tau Z_i - \tau^2/2}]$,
which $\approx 1$ if $Z_i$ is large and positive and $\approx 0$ if $Z_i$ is large and negative.  Seemingly, the posterior distribution contains much more information on inference than the prior distribution does.
Note that this also suggests the (one-sided) clipping hard thresholding, similar to that suggested by Fisher's LDA.
\end{remark}

\subsection{Phase diagram for classification}
We introduce the ARW model for classification:
\[
\sqrt{n}\mu_i   \stackrel{iid}{\sim}    (1-\eps)  \nu_0+\eps  \nu_{\tau},  \qquad 1 \leq i \leq p.
\]
Fix $(\vartheta,r,\theta)$ such that $r>0$, $0< \theta<1$ and $0 < \vartheta<1 - \theta$. Similarly, we  let
\[
\eps=\eps_p=p^{-\vartheta}, \qquad \tau=\tau_p = \sqrt{2r\log(p)} \qquad \mbox{and}\quad n=n_p = p^{\theta}.
\]
It was noted in \cite{JinPNAS, FJY} that
for any fixed $\theta \in (0,1)$, the most  interesting range for $\vartheta$ is $0 < \vartheta < (1 - \theta)$.
When $\vartheta > (1 - \theta)$,   for successful classification, we need $\tau_p \gg \sqrt{\log(p)}$, but this corresponds to Rare/Strong regime, which is relatively easy, for we can separate the nonzero entries of $\mu$ from zero ones by simple thresholding.
For $\rho^*(\cdot)$ be as  in \eqref{phasefunction}, let
\[
\rho^*_{\theta}(\vartheta) =  (1 - \theta)  \rho^*(\frac{\vartheta}{1 - \theta}), \qquad 0 < \vartheta < (1 - \theta),
\]
The following theorem is proved in Fan {\it et al} (2013) \cite{FJY}.
\begin{thm}  \label{thm:class}
Fixing $(\vartheta, \theta, r) \in (0,1)^3$ such that $0 < \vartheta < (1 - \theta)$, and suppose
$\Omega = \Sigma^{-1}$ satisfies (\ref{diagOmega})-(\ref{degreeOmega}), and that the spectral norm of $\Sigma$ is bounded by a constant $C>0$.
  If $r > \rho_{\theta}^*(\vartheta)$, then  the  classification error of the trained HCT classification rule in (\ref{HCTrule}) tends to $0$ as $p \goto \infty$.
If $0 < r < \rho_{\theta}^*(\vartheta)$, then the classification error of any trained classification rule is no less than $1/2 + o(1)$, where $o(1) \goto 0$ as $p \goto \infty$.
\end{thm}

\begin{figure}[tb]
\centering
\includegraphics[width = 2.8 in]{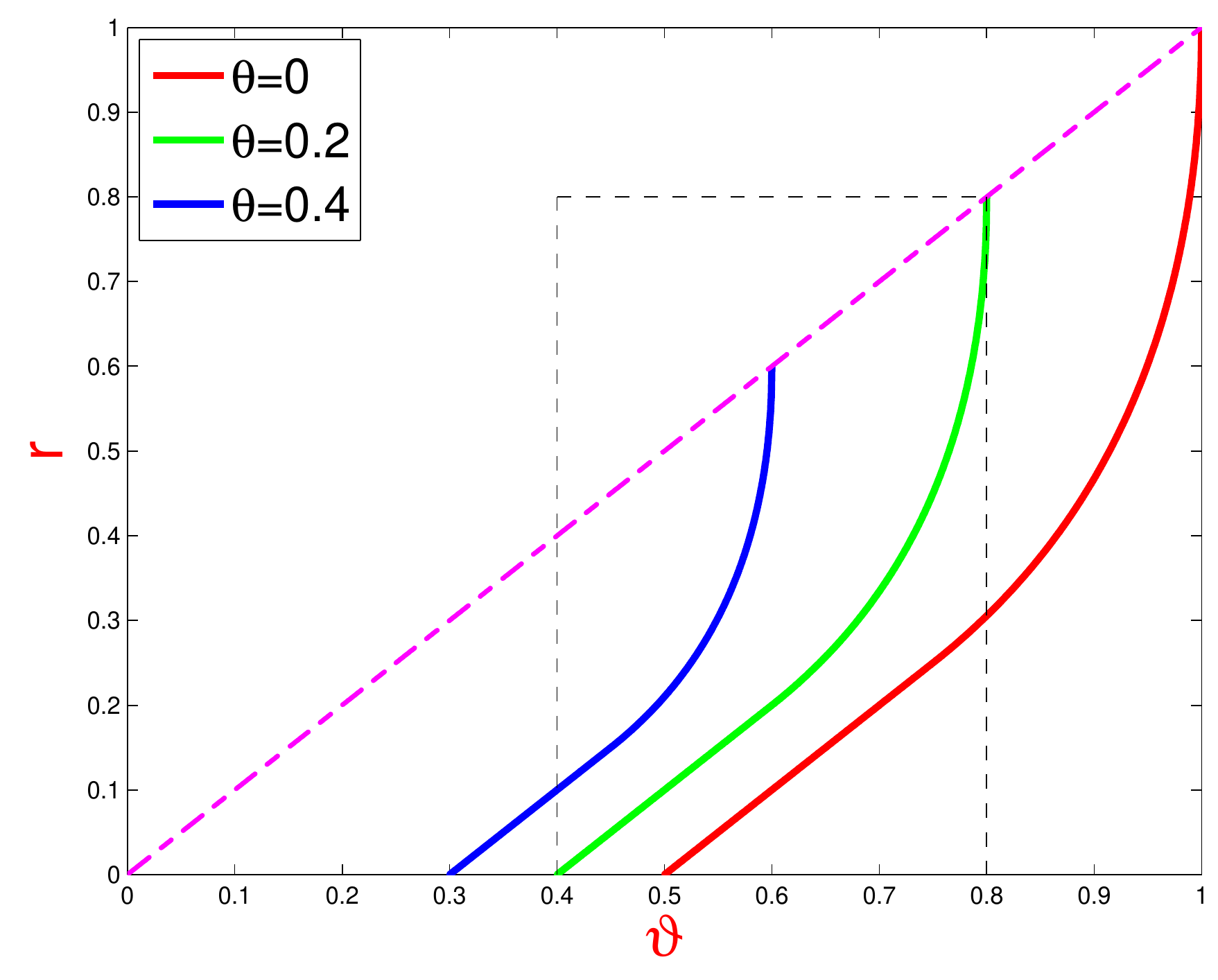}
\includegraphics[width = 2.8 in]{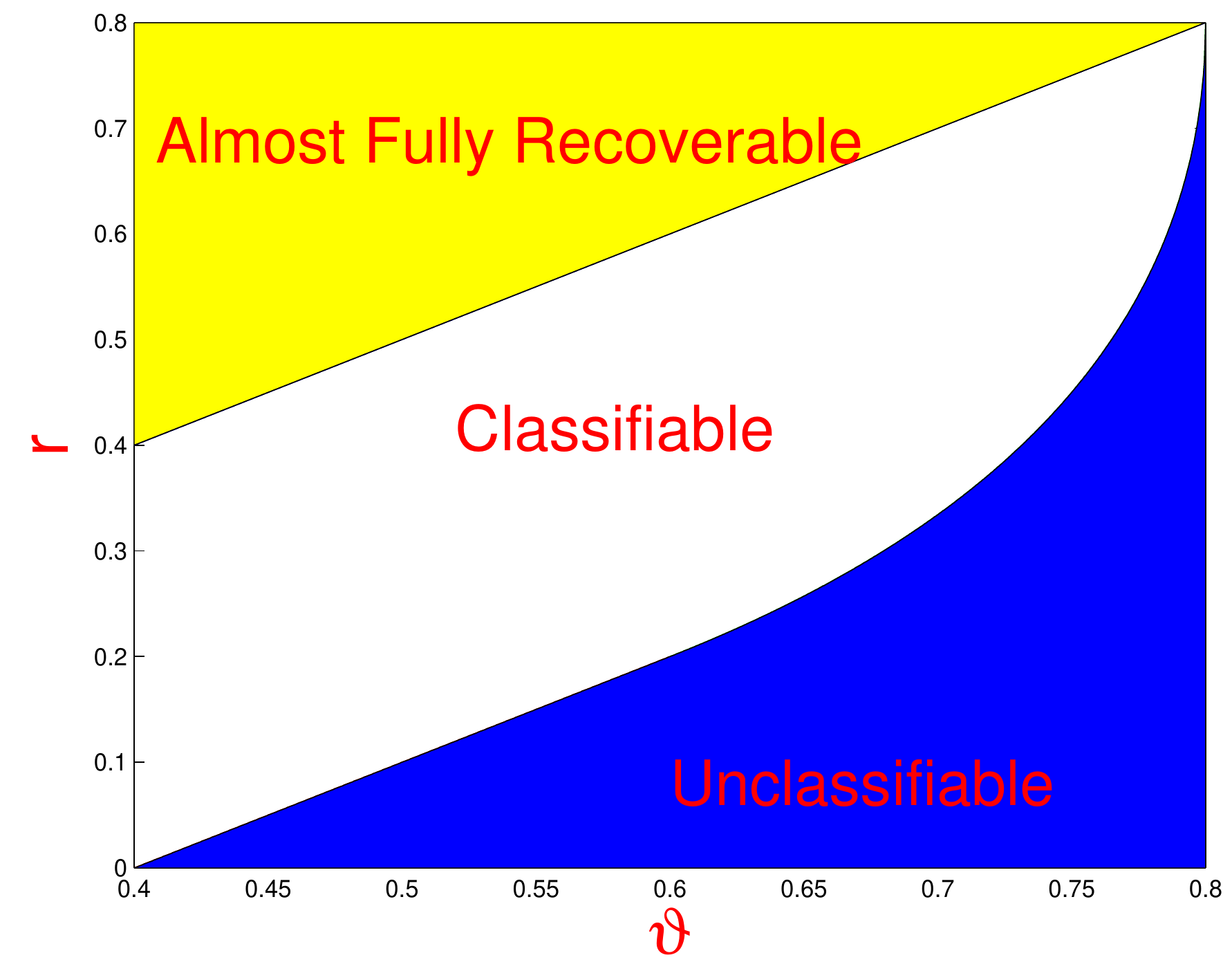}
\caption{Left: curves $r = \rho_{\theta}^*(\vartheta)$ for $\theta = 0, 0.2, 0.4$.
Right: enlargement  the rectangular region bounded by dashed lines  in the left panel.
It partitions into three regions, where in the yellow region, it is not only possible to have successful classifications, but is also possible to separate useful features from useless ones.} \label{fig:Classification}
\end{figure}

There is a similar phase diagram associated with the classification problem.
\begin{itemize}
\item  {\it Region of Classifiable}:  $\{(\vartheta, r):  0 < \vartheta < (1 - \theta),  r > \rho_{\theta}^* (\vartheta)\}$. In this region,  the HC threshold $\hat{t}_p^{HC}$ satisfies $\hat{t}_p^{HC}/ t_p^{ideal} \goto 1$ in probability, where  $t_p^{ideal}$ is the ideal threshold that one would choose if the underlying parameters $(\vartheta, r, \Omega)$ are known.  Also, the classification error of HCT-trained classification rule tends to $0$ as $p \goto \infty$.
\item {\it Region of Unclassifiable}: $\{(\vartheta, r):  0 < \vartheta < (1 - \theta),  r < \rho_{\theta}^* (\vartheta)\}$. In this region,  the classification error of any trained classification rule can not be substantially smaller than $1/2$.
\end{itemize}
See more discussion in \cite{DJ08, DJ09, JinPNAS}. Ingster {\it et al.} (2009)  \cite{IPT} derived independently the classification boundary in a broader setting, but they didn't discuss HC. In Figure~\ref{fig:Classification}, we plot the phase diagrams for $\theta=0, .2, .4$.

The advantage of HC is its optimality in the ARW model.  Note that HCT is a data-driven non-parametric statistic the use of which does not require the knowledge of the ARW parameters.  HC is not tied to the idealized model we discussed here, and can be useful for more general settings. See \cite{DJ08} for applications of HC to cancer classification with
microarray data sets.

Our proposal of threshold choice by HC is very different from  Benjamini-Hochberg's
FDR-controlling method (or Efron's local FDR approach),
where the philosophy is to control the feature FDR (i.e.,  the expected fraction of falsely selected features out of all selected features)   by a small number (e.g., $5\%$). However, this is not necessarily the right strategy when signals are Rare/Weak.
Donoho and Jin (2009) \cite{DJ09} identified a sub-region of Region of Classifiable where to obtain optimal classification behavior,  we must set the feature selection threshold very low
so that  we  include most of the
useful features; but when we do this, we must include
many useless features and the feature FDR is approximately $1$.


\begin{thebibliography}{}
 
\bibitem{Lugosi}
{\sc Addario-Berry},  L.,  {\sc Broutin},   N.,     {\sc Devroye},  L. and {\sc Lugosi},  G. (2010).   On combinatorial testing problems. {\it Ann. Statist}. {\bf 38} 3063--3092.

 
\bibitem{Erycluster}
{\sc Arias-Castro},  E.,  {\sc Candes},  E.  and {\sc Durand}, A.   (2011).
Detection of an anomalous cluster in a network
{\it Ann. Statist}. {\bf 39}(1)   278--304.

 
\bibitem{Castro2}
{\sc Arias-Castro},  E.,  {\sc Candes},  E.  and {\sc Plan},  Y.  (2011).  Global testing under sparse alternatives: ANOVA, multiple comparisons and the higher criticism
{\it Ann. Statist}. {\bf 39}(5) 2533--2556.


\bibitem{WangMeng}
 {\sc Arias-Castro},  E.  and {\sc Wang},    M.  (2013).   Distribution free
 tests for sparse heterogeneous mixtures. {\it arXiv:1308.0346}.

 

\bibitem{BH95}
{\sc Benjamini},  Y. and {\sc Hochberg},   Y.   (1995).     Controlling the false discovery rate: A
practical and powerful approach to multiple testing. {\it J. Roy. Statist. Soc. B}
{\bf 57}  289--300.

\bibitem{Bennett}
{\sc Bennett},  M.F.,  {\sc Melatos},  A.,  {\sc Delaigle},  A. and {\sc Hall}.  P.  (2012).  Reanalysis of $F$-statistics
gravitational-wave search with the higher criticism statistics. {\it Astrophys. J.}  {\bf 766}(99)  1-10.

 

\bibitem{Bogdan1}
{\sc Bogdan},   M.,    {\sc Chakrabarti},  A.,  {\sc Frommlet},  F. and  {\sc Ghosh}, J.  (2011).
Asymptotic Bayes-optimality under sparsity of some multiple testing procedures.
{\it Ann. Statist.} {\bf 39}(3)  1551--1579.

\bibitem{Bogdan2}
{\sc Bogdan},  M.,    {\sc Ghosh}, J. and  {\sc Tokdar},  T.   (2008).
A comparison of the Benjamini-Hochberg procedure with some Bayesian rules for
multiple testing. {\it Beyond Parametrics in Interdisciplinary Research: Festschrift in Honor
of Professor Pranab K. Sen}  {\bf 1}  211--230.


\bibitem{Box}
{\sc Box}, M.  and {\sc Meyer}, D. (1986).
An analysis for unreplicated fractional factorials. {\it  Technometrics} {\bf 28}, 11--18.

\bibitem{CJJ}
{\sc Cai},  T., {\sc  Jeng},  J. and {\sc Jin},  J.   (2011).
Detecting sparse heterogeneous and heteroscedastic mixtures.
 {\it J. Roy. Statist. Soc. B}.  {\bf 73}  629--662.


\bibitem{CJL}
{\sc Cai},   T.,    {\sc  Jin},  J. and {\sc  Low}, M.  (2007).
Estimation and confidence sets for sparse normal mixtures.
{\it Ann. Statist}.  {\bf 35}(6)  2421--2449.


\bibitem{CaiLuo}
{\sc Cai},  T., {\sc  Liu},  W. and {\sc Luo},  X. (2010).
A constrained $L^1$ minimization approach to sparse precision matrix estimation.
{\it J. Amer. Statist. Assoc}. {\bf 106} 594--607.


\bibitem{CaiWu}
{\sc Cai},  T. and {\sc Wu},  Y.  (2012).  Optimal detection
for sparse mixtures.  {\it arXiv:1211.2265}.



\bibitem{Cayon2}
{\sc Cayon}, L. and {\sc Banday},  A.J.  {\it et al}. (2006).
No Higher Criticism of the Bianchi-corrected Wilkinson Microwave Anisotropy Probe data.
{\it Mon. Not. Roy. Astron. Soc.}  {\bf 369}(2)   598--602.

\bibitem{Cayon1}
{\sc  Cayon},  L., {\sc  Jin},  J. and  {\sc Treaster},  A.  (2004).   Higher Criticism statistic: detecting and identifying non-Gaussianity in the WMAP first year data. {\it Mon. Not. Roy. Astron. Soc}.  {\bf 362}  826--832.

 

\bibitem{Cruz}
{\sc Cruz}, M.,    {\sc Cayon},  L.,     {\sc Martinez-Gonzalez},  E.,    {\sc Vielva},  P.  and {\sc Jin},  J.
(2007).   The non-Gaussian cold spot in the 3 year Wilkinson Microwave Anisotropy Probe data.
{\it  Astrophys. J.}  {\bf 655}(1)  11--20.

 
 \bibitem{delaCruz}
{\sc De la Cruz}, O., {\sc Wen}, X.,  {\sc Ke}, B.,  {\sc Song}, M. and {\sc Nicolae}, D. L. (2010).   Gene, region and pathway level analyses in whole-genome studies. {\it Genet. Epidemiol}. {\bf 34}  222-231.



\bibitem{Delaigle}
{\sc Delaigle},  A. and  {\sc Hall},  P.  (2009).
{\it Higher criticism in the context of unknown distribution, non-independence and classification}. In Perspectives in Math- ematical Sciences I: Probability and Statistics, 109--138 (eds N. Sastry, M. Delampady, B. Rajeev and T.S.S.R.K. Rao. World Scientific.

\bibitem{Robustness}
{\sc Delaigle},   A.,   {\sc  Hall},  J. and {\sc Jin},  J. (2011).
Robustness and accuracy of methods for high dimensional data analysis based on Student's t statistic. {\it J. Roy. Statist. Soc. B}.   {\bf 73}  283--301.


 
\bibitem{Una}
{\sc de Una-Alvarez},   J.  (2012).   The Beta-Binomial SGoF method for multiple dependent tests.
{\it Statistical applications in genetics and molecular biology}.   {\bf 11}(3)  1544-6115.


\bibitem{Nissim}
{\sc Dinur}, I. and {\sc Nissim},   K.   (2003). Revealing information while preserving privacy. {\it Proceedings of the twenty-second ACM SIGMOD-SIGACT-SIGART symposium on Principles of database systems} 202--210  ACM Press.

\bibitem{Donoho}
{\sc Donoho},  D.  (2006).   Compressed sensing. {\it IEEE Trans. Inform. Theory}  {\bf 52}(4)
1289--1306.

\bibitem{DJ04}
{\sc Donoho},   D. and {\sc Jin},  J. (2004).  Higher criticism for detecting
sparse heterogeneous mixtures.  {\it Ann. Statist}.  {\bf 32}   962-994.


\bibitem{DJ08}
{\sc Donoho},  D. and {\sc Jin},  J.  (2008).  Higher Criticism thresholding: optimal feature selection when useful features and rare and weak.
{\it Proc. Natl. Acad. Sci.} {\bf 105}(39)   14790--14795.

\bibitem{DJ09}
{\sc Donoho},   D. and {\sc Jin},   J.  (2009).
Feature selection by Higher Criticism thresholding: optimal phase diagram. {\it Phil. Tran. Roy. Soc. A} {\bf 367}
4449--4470.


\bibitem{DonoJohn94}
{\sc Donoho}, D. and  {\sc Johnstone}, I. (1994).
Ideal spatial adaptation by wavelet shrinkage. {\it Biometrika} {\bf 81}(3) 425--455.

\bibitem{DonoJohn95}
{\sc Donoho}, D. and  {\sc Johnstone}, I. (1994).
Adapting to unknown smoothness via wavelet shrinkage. {\it J. Amer. Statist. Assoc.} {\bf 90}(432) 1200--1224.



\bibitem{DMM09}
{\sc Donoho}, D., {\sc Maleki}, A. and {\sc Montanari}, A. (2009).
Message-passing algorithms for compressed sensing. {\it Proc. Natl.  Acad. Sci.} {\bf 106}(45) 18914--18919.

\bibitem{DonohoStark}
{\sc Donoho}, D., {\sc Stark}, P.  (1989).
Uncertainty principles and signal recovery. {\it SIAM J. Appl. Math.} {\bf 49}(3) 906--931.

 

\bibitem{EfronNull}
{\sc Efron},  B.  (2004).   Large-scale simultaneous hypothesis testing:
the choice of a null hypothesis.  {\it J. Amer. Statist. Assoc.}  {\bf 99}  96--104.

\bibitem{EfronLSI}
{\sc Efron},  B.  (2011).   {\it Large-Scale Inference: Empirical Bayes methods for
estimation, testing, and prediction}.  IMS Monographs, Cambridge Press.


\bibitem{LARS}
{\sc Efron}, B., {\sc Hastie}, T., {\sc Johnstone}, I. and {\sc Tibshirani}, R. (2004)
Least angle regression. {\it Ann. Statist.} {\bf 32} (2) 407--840.

 
\bibitem{FanLi}
{\sc Fan}, J. and {\sc Li}, R. (2001).
Variable selection via nonconcave penalized likelihood and its oracle properties. {\it J. Amer. Statist. Assoc.} {\bf 96}(456) 1348--1360.

\bibitem{FanLv}
{\sc Fan}, J. and {\sc Lv}, J. (2008). Sure independence screening for ultrahigh dimensional feature space.
{\it J. Roy. Statist. Soc. B.} {\bf 70} (5) 849--911.


\bibitem{ISIS}
{\sc Fan}, J., {\sc Samworth}, R. and {\sc Wu}, Y. (2009).
Ultrahigh dimensional feature selection: beyond the linear model.
{\it J. Mach. Learn. Res.} {\bf 10} 2013--2038.

\bibitem{FJY}
{\sc Fan},  Y.,  {\sc Jin}, J.  and {\sc Yao},  Z.  (2013).
Optimal classification by Higher Criticism in sparse Gaussian graphic  model.
{\it Ann. Statist.}, {\bf 41}(5), 2263--2702.

\bibitem{Fienberg}
{\sc Fienberg},  S. and {\sc Jin},  J.  (2012).   Privacy-preserving data sharing in high dimensional regression and classification settings. {\it J. Privacy and Confidentiality}  {\bf 4}(1)  Article 10.

\bibitem{glasso}
{\sc Friedman},  J., {\sc   Hastie},  T. and {\sc Tibshirani},   R.  (2007). Sparse inverse covariance
estimation with the graphical lasso. {\it Biostatistics} {\bf 9} 432-441.

\bibitem{Frieze}
{\sc Frieze}, A.M. and  {\sc Molloy}, M. (1999).   Splitting an expander graph. {\it J.   Algorithms}, {\bf 33}(1),   166--172.

\bibitem{Ingstervariable}
{\sc Gayraud}, G. and  {Ingster}, Y.I.  (2011).
Detection of sparse variable functions. {\it  arXiv:1011.6369}.


\bibitem{Ge}
{\sc Ge},  Y.  and {\sc Li},  X.  (2012).
Control of the False Discovery Proportion for
independently tested null hypotheses. {\it J. Probab. and Statist.}
  2012   Article ID 320425  19 pages.


\bibitem{GJW}
{\sc Genovese},  C.,  {\sc  Jin},  J.,  {\sc  Wasserman},  L. and {\sc Yao},  Z.   (2012).
A comparison of the lasso and marginal regression.
{\it J. Mach. Learn. Res.}  {\bf 13}  2107--2143.


 
\bibitem{Greenshtein}
{\sc Greenshtein},   E.  and {\sc Park},   J.   (2012).  Robust test for detecting a signal in a high dimensional
sparse normal vector. {\it J. Statist. Planning and Inference}   {\bf 142}  1445--1456.



\bibitem{HJ08}
{\sc Hall},   P.   and  {\sc Jin},  J.   (2008).    Properties of higher criticism under strong dependence.
{\it Ann. Statist}. {\bf 36}(1) 381--402.


\bibitem{HJ09}
{\sc Hall},  P. and {\sc Jin},   J.  (2010).  Innovated Higher Criticism for detecting
sparse signals in correlated noise.  {\it Ann. Statist.}
{\bf 38}(3) 1686--1732.

 
\bibitem{Haupt1}
{\sc Haupt},  J., {\sc  Castro},   R. and {\sc Nowak}, R.  (2008).
Adaptive discovery of sparse signals in noise.
{\it Signals, Systems and Computers, 2008 42nd Asilomar Conference}.  1727--1731.

\bibitem{Haupt2}
{\sc Haupt}, J., {\sc  Castro},  R.  and {\sc Nowak},  R. (2010).
 Improved bounds for sparse recovery from adaptive measurements.
 {\it Information Theory Proceedings (ISIT), 2010} 1565--1567.

 

\bibitem{GWAS2}
{\sc He},  S. and {\sc Wu},  Z.   (2011).
Gene-based Higher Criticism methods for large-scale exonic single-nucleotide
polymorphism data. {\it BMC Proceedings} {\it 5} (Suppl 9):S65.

 
\bibitem{Ingster97}
{\sc Ingster}, Y.I. (1997). Some problems of hypothesis testing leading to infinitely divisible distribution. {\it   Math. Methods Statist.}   {\bf 6}  47--69.

\bibitem{Ingster99}
{\sc Ingster},  Y.I.  (1999).  Minimax detection of a signal for $l^p_n$-balls.  {\it  Math. Methods Statist.}  {\bf 7} 401--428.

\bibitem{IPT}
{\sc Ingster},  Y.I.,  {\sc Pouet},  C. and {\sc Tsybakov},  A.  (2009).
Classification of sparse high-dimensional vectors.  {\it Phil. Trans. R. Soc. A}
{\bf 367} 4427--4448.

\bibitem{IPV}
{\sc Ingster},  Y.I.,  {\sc Tsybakov}, A. and  {\sc Verzelen},  N.  (2010).
Detection boundary in sparse regression.
{\it Electron. J. Statist}. {\bf 4}  1476--1526.


\bibitem{Wellner2004}
{\sc Jager},  L. and  {\sc Wellner},  J.A. (2004).
On the ``Poisson boundaries" of the family of weighted Kolmogorov statistics.
{\it IMS Monograph}  {\bf 45}  319--331.


\bibitem{Wellner2007}
{\sc Jager},  L. and {\sc Wellner},  J.A.  (2007).
Goodness-of-fit tests via phi-divergences. {\it Ann. Statist.}  {\bf 35}(5)
2018--2053.

 
\bibitem{Jeng1}
{\sc Jeng},  J.,  {\sc Cai},  T. and {\sc Li},  H.  (2010).
Optimal sparse segment identification with application in copy number variation analysis.
{\it J. Amer. Statist. Assoc.} {\bf 105}(491)  1156--1166.

\bibitem{Jeng2}
{\sc Jeng},   J.,    {\sc Cai},  T. and {\sc Li},  H.   (2013).
Simultaneous discovery of rare and common segment variants.
{\it Biometrika} {\bf 100}(1) 157--172.



\bibitem{JiJin}
{\sc Ji}, P. and {\sc Jin}, J. (2011). UPS delivers optimal phase diagram in high dimensional variable selection.
{\it  Ann. Statist.}  {\bf 40}(1)   73-103.

\bibitem{JinThesis}
{\sc Jin},  J.    (2003).    Detecting and estimating sparse mixtures.   {\it Ph.D. Thesis},  Department of Statistics, Stanford University.

 
\bibitem{JASA}
{\sc Jin}, J. (2012). Comment on ``Estimating false discovery proportion under arbitrary
covariance dependence". {\it J. Amer. Statist. Assoc.} {\bf 107} 1042--1045.


 
\bibitem{JinPNAS}
{\sc Jin},  J. (2009).   Impossibility of successful classification when useful features are rare and weak.
{\it Proc. Natl.  Acad. Sci.} {\bf 106}(22) 8859--9964.

 

\bibitem{Starck}
{\sc Jin},    J.,      {\sc  Starck},  J.-L.,     {\sc Donoho},   D., {\sc  Aghanim},   N. and  {\sc Forni},   O.  (2005).   Cosmological non-gaussian signature detection: Comparing performance of different statistical tests.   {\it EURASIP J. Appl. Signal Processing} {\bf 15}  2470--2485.


\bibitem{JW}
{\sc Jin},  J. and {\sc Wang},  W.  (2012) Optimal feature selection by Higher Criticism thresholding in
spectral clustering. {\it Manuscript}.

\bibitem{JZZ}
{\sc Jin}, J., {\sc Zhang}, C.-H. and {\sc Zhang}, Q. (2012). Optimality of Graphlet Screening in high dimensional variable selection.  {\it arXiv:1204.6452}.

 

\bibitem{Ke}
{\sc Ke},  Z.,    {\sc  Jin},  J. and {\sc Fan},  J.  (2012).
Covariance assisted screening and estimation. {\it arXiv:1205.4645}.


\bibitem{Kendall}
{\sc Kendall},   D.    (1980).   Discussion of  ``simulating the Ley hunter" by Simon Broadbent  {\it J. Royal Stat. Soc. A}   {\bf 2} 109--140.


 

\bibitem{Laurent}
{\sc Laurent},  B.,  {\sc  Marteau}, C.  and {\sc Maugis-Rabusseau},  C. (2013).
Non-asymptotic detection of two-component mixture with unknown means. {\it aiXiv:1304.6924}.

 

\bibitem{ShaoQM}
{\sc Liu},  W.  and {\sc Shao}.  Q.M. (2013).  A Cram\'er Rao moderate deviation theorem for
Hotelling's $T^2$-statistic with applications to global tests.   {\it Ann. Statist.} {\bf 41}(1)  296--322.


\bibitem{Martin}
{\sc Martin},  L., {\sc Gao},  G.,  {\sc  Kang},  G.,  {\sc Fang},  Y. and {\sc  Woo}, J.   (2009).
Improving the signal-to-noise ratio in genome-wide association studies.
{\it Genetic Epidemiology}  {\bf 33} Suppl 1 29--32.


\bibitem{Meinshausen}
{\sc Meinshausen},  N.  (2008).
Hierarchical testing of variable importance.
{\it  Biometrika}  {\bf 95}(2)   265--278.

\bibitem{Rice}
{\sc Meinshausen},  N. and {\sc Rice},  J.  (2006).
Estimating the proportion of false null hypotheses among a large number of independently tested hypotheses.
{\it Ann. Statist.} {\bf 34}(1)  373--393.


\bibitem{Lin}
{\sc Mukherjee},  R.,   {\sc Pillai},  N.   and {\sc Lin}, X.  (2013).
Hypothesis testing for sparse binary regression.  {\it arXiv:1308.0764}.




\bibitem{Neill}
{\sc Neill},  D.  (2006).
Detection of spatial and spatio-temporal Clusters.
{\it Ph.D Thesis}, School of Computer Science,  Carnegie Mellon University.


\bibitem{Neill0}
{\sc Neill},  D. and  {\sc Lingwall},  J.    (2007).
A nonparametric scan statistic for multivariate disease surveillance.
{\it Advances in Disease Surveillance}   {\bf 4} 106-106.

\bibitem{Ghosh}
{\sc Park},  J. and {\sc Ghosh},  J. (2010).
A guided random walk through some high dimensional problems.
{\it Sankhya}  {\bf 72-A} 81--100.

\bibitem{GWAS1}
{\sc Parkhomenko},   E.,  {\sc Tritchler},  D.  and {\sc Lemire},  M.  {\it et al}.  (2009).
Using a higher criticism statistic to detect modest effects in a genome-wide study of rheumatoid arthritis.
{\it BMC Proceedings}  {\bf 3} (Suppl 7):S40.

\bibitem{GRN}
{\sc Peng}, J.,  {\sc Wang}, P.,  {\sc Zhou}, N. and {\sc Zhu}, J. (2010). Partial correlation estimation by joint sparse regression model. {\it  J. Amer. Statist. Assoc}., {\bf 104}(486),735--746.



\bibitem{Pires}
{\sc Pires},  S.,    {\sc Starck},  J.-L.,   {\sc Amara}, A.,  {\sc Refregier}, A. and {\sc Teyssier},  R. (2009).
Cosmological models discrimination with Weak Lensing. {\it Astron. Astrophys.}   {\bf 505}  969--979.


\bibitem{RW}
{\sc Roeder}, K. and  {\sc Wasserman},  L.  (2009).
Genome-wide significance levels and weighted hypothesis testing.  {\it Statist. Sci.}  {\bf 24}(4)  398--413.



\bibitem{Sabatti}
{\sc Sabatti},  C.,  {\sc Service},  S. and {\sc Hartikainen},  A. L.   {\it et al.} (2008).
Genome-wide association analysis of metabolic traits in a birth cohort from a founder population.
{\it Nature Genetics}  {\bf  41} 35--46.

\bibitem{anomaly}
{\sc Saligrama},  V. and {\sc Zhao},   M.  (2012).  Local anomaly detection. {\it AISTATS 2012}.

 
\bibitem{Tib96}
{\sc Tibshirani}, R. (1996).
Regression shrinkage and selection via the lasso.
{\it J. Roy. Statist. Soc. B.} {\bf 58} 267--288.

\bibitem{Tibs}
{\sc Tibshirani},    R.,  {\sc  Hastie},   T.,  {\sc  Narasimhan},   B. and {\sc  Chu},   G.  (2002).   Diagnosis of multiple cancer types by shrunken centroids of gene expression.    {\it Proc. Natl. Acad. Sci.}  {\bf 99}   6567--6572.


\bibitem{Tukey76}
{\sc Tukey},  J.W. (1976).  T13 N: The Higher Criticism.   {\it  Course notes, Stat 411, 1976. Princeton University}.


\bibitem{Tukey89}
{\sc Tukey},  J.W. (1989).   Higher Criticism for individual significances in several tables or parts of tables. {\it Internal working paper, Princeton University}.


\bibitem{Tukey94}
{\sc Tukey},   J.W. (1994).  {\it The problem of Multiple Comparisons in The Collected Works of John W. Tukey}.  Vol.  III, 1948--1983, Edited by Henry I. Braun. (Original manuscript 1953.) Chapman \& Hall.

 

\bibitem{Vielva}
{\sc Vielva}, P.   (2010).
A comprehensive overview of the cold spot.
{\it Advances in Astronomy}.  2010  Article ID 592094  20 pages.


 
\bibitem{Wassermanbook}
{\sc Wasserman}, L. (2006). {\it All of nonparametric statistics}. Springer, NY.


\bibitem{WellnerKoltchinskii}
{\sc Weller},  J.A. and {\sc Koltchinskii},  V.  (2003).   A note on the
asymptotic distribution of Berk-Jones type statistics under the null
hypothesis.  In {\it High Dimensional Probability III} (J. Hoffmann-Jorgensen, M. B. Marcus and J.A. Wellner, eds.) Birkh\"auser, Basel.

 
\bibitem{Wu}
{\sc Wu},   Z.,     {\sc Sun},  Y.,   {\sc  He},  S.,  {\sc Choy},  J.,  {\sc Zhao},  H.   and {\sc Jin},  J.  (2012).
Detection boundary and Higher Criticism approach for
sparse and weak genetic effects.   {\it Manuscript}.

 
\bibitem{MCP}
{\sc Zhang}, C.-H. (2010).
Nearly unbiased variable selection under minimax concave penalty.
{\it Ann. Statist.},  {\bf 38} (2),  894--942.



\bibitem{FoBo}
{\sc Zhang}, T. (2011) Adaptive forward-backward greedy algorithm for learning sparse representations.
{\it IEEE Trans. Inform. Theory} {\bf 57}(7) 4689--4708.


\bibitem{Chen}
{\sc Zhong},  P.,  {\sc  Chen},  S.  and {\sc  Xu},   M.  (2013).
Test alternative to higher criticism for high dimensional means under sparsity and column-wise dependence.
{\it Ann. Statist.}, to appear.

\end{thebibliography}
\end{document}